\documentclass{article}

\usepackage{microtype}
\usepackage{graphicx}
\usepackage{subfigure}
\usepackage{booktabs} 

\usepackage{algorithm, algorithmic}
\usepackage{tikz} 
\usepackage{scrextend}
\usepackage{amsmath,amsthm,amssymb,url,amscd}
\usepackage{amsfonts}
\usepackage{multirow}
\usepackage{enumerate}
\usepackage{accents}
\usepackage{amscd}
\usepackage{subfigure}
\usepackage{float}
\usepackage{bm}
\usepackage{caption}
\usepackage[colorlinks=true, linkcolor = blue, citecolor=tumred]{hyperref}
\usepackage{cleveref}
\usepackage{xcolor}
\usepackage{setspace}
\usepackage{pgfplots}
\pgfplotsset{compat=newest}
\usetikzlibrary{plotmarks}
\usetikzlibrary{arrows.meta}
\usepgfplotslibrary{patchplots}
  
\pgfplotsset{plot coordinates/math parser=false}
\newlength\figureheight
\newlength\figurewidth

\DeclareMathOperator{\Id}{Id}

\DeclareMathOperator*{\argmin}{arg\,min}
\DeclareMathOperator{\vecc}{vec}
\DeclareMathOperator{\rank}{rank}

\DeclareMathOperator{\R}{\mathbb{R}}

\DeclareMathOperator{\mat}{mat}

\DeclareMathOperator{\diag}{diag}
\DeclareMathOperator{\dg}{dg}

\newcommand\f[1]{\mathbf{#1}} 

\newcommand\Rdd{\R^{d_1 \times d_2}}

\newcommand\hk{^{(k)}}
\newcommand\hkk{^{(k+1)}}

\newtheorem{theorem}{Theorem}[section]

\usepackage{hyperref}

 \usepackage[accepted]{icml2020}

\icmltitlerunning{Scalable Second-Order Optimization for Ill-Conditioned Matrix Completion}

\begin{document}

\twocolumn[
\icmltitle{Escaping Saddle Points in Ill-Conditioned Matrix Completion\\ with a Scalable Second Order Method}
\icmlsetsymbol{equal}{*}
 
\begin{icmlauthorlist}
\icmlauthor{Christian K\"ummerle}{jhu}
\icmlauthor{Claudio M. Verdun}{tum}
\end{icmlauthorlist}

\icmlaffiliation{jhu}{Department of Applied Mathematics \& Statistics, Johns Hopkins University, Baltimore, USA}
\icmlaffiliation{tum}{Department of Mathematics and Department of Electrical and Computer Engineering, Technical University of Munich, Munich, Germany, e-mail: \texttt{verdun@ma.tum.de}}

\icmlcorrespondingauthor{Christian K\"ummerle}{\texttt{kuemmerle@jhu.edu}}

\icmlkeywords{Low-Rank Matrix Completion, Low-Rank Matrix Recovery, Compressive Sensing, Iteratively Reweighted Least Squares, Saddle Points, Second Order Methods, Ill-conditioned matrices, Majorization Minimization}

\vskip 0.3in
]



\printAffiliationsAndNotice{}  

\begin{abstract}

    We propose an iterative algorithm for low-rank matrix completion that can be interpreted as both an iteratively reweighted least squares (IRLS) algorithm and a saddle-escaping smoothing Newton method applied to a non-convex rank surrogate objective. It combines the favorable data efficiency of previous IRLS approaches with an improved scalability by several orders of magnitude. Our method attains a local quadratic convergence rate already for a number of samples that is close to the information theoretical limit. We show in numerical experiments that unlike many state-of-the-art approaches, our approach is able to complete very ill-conditioned matrices with a condition number of up to $10^{10}$ from few samples.

\end{abstract}

\vspace{-8mm}
\section{Introduction}\label{intro}
\vspace{-2mm}

In different areas of machine learning and signal \mbox{processing}, low-rank models have turned out to be a powerful tool for the acquisition, storage and computation of \mbox{information}. 
In many of these applications, an important sub-problem is to infer the low-rank model from partial or incomplete data \cite{Davenport16,ChiLuChen19}. 

This problem is called \emph{low-rank matrix completion}: Given a matrix $\f{X}^0 \in \Rdd$ of rank-$r$ and an index set $\Omega \subset [d_1] \times [d_2]$, the task is to reconstruct $\f{X}^0$ just from the knowledge of $\Omega$ and $P_{\Omega}(\f{X}^0)$, where $P_{\Omega}: \Rdd \to \mathbb{R}^m$ is the subsampling operator that maps a matrix to the entries indexed by $\Omega$. It is well-known that this can be reformulated as the NP-hard \emph{rank minimization} problem \cite{Recht10}
\vspace{-2mm}
\begin{equation}\label{rank_equation}
    \min_{\f{X} \in \R^{d_1 \times d_2}} \rank(\f{X}) \quad \mbox{ subject to } P_{\Omega}(\f{X}) = P_{\Omega}(\f{X}^0).
\end{equation}

From an optimization point of view, \cref{rank_equation} is particularly difficult to handle due to two properties: its \emph{non-convexity} and its \emph{non-smoothness}. A widely studied approach in the literature is to replace the $\rank(\f{X})$ by the convex nuclear norm $\|\f{X}\|_* = \sum_{i=1}^d \sigma_i(\f{X})$ \cite{FHB03}. However, such a convex relaxation approach has two main drawbacks: it is \emph{computational demanding} for large problems, since it is equivalent to a semidefinite program \cite{Recht10}, and it is not \emph{data efficient}, since the nuclear norm minimizer under the affine constraint might not coincide with the minimizer $\f{X}^0$ of \cref{rank_equation} if the number of samples $m$ is just slightly larger than the number of the degrees of freedom of $\f{X}_0$ \cite{Amelunxen14}.

To overcome these drawbacks, a variety of alternative approaches have been proposed and studied, many of which optimize an empirical loss defined for a matrix factorization model, or which optimize the empirical loss using Riemannian manifold structures, see \cite{ChiLuChen19} for a recent survey. These approaches are much more scalable, and furthermore, are often able to reconstruct the low-rank matrix from fewer samples than a convex formulation.

However, a closer inspection of the theoretical guarantees of these algorithms suggests that the performance of those algorithms deteriorates as the \emph{condition number} $\kappa$ of $\f{X}^0$ increases. For example, if $D=\max(d_1,d_2)$, the sufficient condition on the required random samples for the Riemannian gradient descent algorithm  \cite{wei_cai_chan_leung} is  $m=\Omega(\kappa^6 D r^2 \log D)$, and a polynomial dependence on $\kappa$ can be found for results on gradient descent for matrix factorization \cite{ChiLuChen19}.\footnote{An exception is the result of \cite{hardt_wotters}, whose sample complexity exhibits a \emph{logarithmic} dependence on $\kappa$. On the other hand, its dependence on the rank $r$ is a high-order polynomial.}

Retrieving ill-conditioned matrices from partial information is a problem that arises in many important areas such as the discretization of PDE-based inverse problems with Fredholm equations \cite{cloninger_czaja_bai_basser} or spectral estimation problems, where Hankel matrices with condition number $\approx 10^{15}$ may appear \cite{liao_fannjiang_music16}. Not only theoretically, but also in practice, non-convex approaches often struggle to recover such ill-conditioned matrices. To overcome this, this paper proposes \texttt{MatrixIRLS}, an efficient second order least squares approach that aims to solve \emph{ill-conditioned} matrix completion problems that are \emph{statistically hard}, i.e., problems in which just very few entries are known. In the following, let $d=\min(d_1,d_2)$.

\vspace{-4mm}
\section{MatrixIRLS for log-det rank surrogate}\label{surrogate}
\vspace{-2mm}
We propose a method that is based on the minimization of quadratic models of smoothed log-det objectives to obtain a scalable, but unbiased method for rank minimization \eqref{rank_equation}. It has already been observed in several papers \cite{Fazel02,Candes13} that optimizing the smoothed log-det objective $\sum_{i=1}^d \sigma_i(\f{X}+\epsilon \f{I})$ for some $\epsilon > 0$ can lead to minimum rank solutions. In particular, it can be shown that the minimizer of the smoothed log-det objective coincides as least as often with the rank minimizer as the convex nuclear norm minimizer \cite{Foucart18}.

However, finding the global minimizer of a non-convex and non-smooth rank surrogate can be very challenging, as the existence of ample sub-optimal local minima might deter the success of many local optimization approaches. Furthermore, applications such as in recommender systems \cite{koren_bell_volinsky} require solving very high-dimensional problem instances so that it is impossible to store full matrices, let alone to calculate many singular values of these matrices.

Let now $\epsilon > 0$ and $F_{\epsilon}:\Rdd \to \R$ be the \emph{smoothed log-det objective} defined as
$F_{\epsilon}(\f{X}):= \sum_{i=1}^d f_{\epsilon}(\sigma_i(\f{X}))$
with 
\vspace*{-2mm}
\begin{equation} \label{eq:smoothing:Fpeps}	
f_{\epsilon}(\sigma) =
\begin{cases}
 \log|\sigma|, & \text{ if } \sigma \geq \epsilon, \\
 \log(\epsilon) + \frac{1}{2}\Big( \frac{\sigma^2}{\epsilon^2}-1\Big), & \text{ if } \sigma < \epsilon.
 \end{cases}
  \vspace*{-1mm}
 \end{equation}
It can be shown that that $F_{\epsilon}$ is continuously differentiable with $\epsilon^{-2}$-Lipschitz gradient \cite{AnderssonCarlssonPerfekt16}. It is clear that the optimization landscape and of $F_{\epsilon}$ crucially depends on the smoothing parameter $\epsilon$. 
 
We propose now an iterative algorithm that minimizes quadratic upper bounds of $F_{\epsilon}$ under the data constraint to provide a descent in $F_{\epsilon}$, before updating the smoothing parameter $\epsilon$. It can be interpreted within the framework of \emph{Iteratively Reweighted Least Squares (IRLS)} \cite{Weiszfeld37,Fornasier11,Mohan10,KS18}, as the main step of each iteration can be regarded as solving a weighted least squares problem.

The precise shape of the quadratic model can be described by a \emph{weight operator} $W\hk$, which we define as follows: Let $\epsilon_k > 0$ and $\f{X}\hk \in \Rdd$ be a matrix with singular value decomposition $\f{X}\hk =  \f{U}_k \dg(\sigma^{(k)}) \f{V}_k^{*}$, i.e., $\f{U}_k \in \R^{d_1 \times d_1}$ and $\f{V}_k \in \R^{d_2 \times d_2}$ are orthonormal matrices. Then we define the linear operator $W\hk: \Rdd \to \Rdd$ such that 
\begin{equation} \label{eq:def:W}
	W^{(k)}(\f{Z}) = \f{U}_k \left[\f{H}_k \circ (\f{U}_k^{*} \f{Z} \f{V}_k)\right] \f{V}_k^{*},
\end{equation}
where $\f{H}_k \circ (\f{U}_k^{*} \f{Z} \f{V}_k)$ is the entrywise product of $\f{H}_k$ and $\f{U}_k^{*} \f{Z} \f{V}_k$, and $\f{H}_k \in \R^{d_1 \times d_2}$ is a matrix with positive entries such that
$(\f{H}_k)_{ij} := \left(\max(\sigma_i^{(k)},\epsilon_k) \max(\sigma_j^{(k)},\epsilon_k)\right)^{-1}$
for all $i \in [d_1]$ and $j \in [d_2]$. The weight operator $W^{(k)}$ is a positive, self-adjoint operator with strictly positive eigenvalues that coincide with the entries of the matrix $\f{H}_k \in \Rdd$. Using the definition of $W\hk$, we describe \texttt{MatrixIRLS} in \Cref{algo:MatrixIRLS}.

\begin{algorithm}[tb]
\caption{\texttt{MatrixIRLS} for low-rank matrix recovery} \label{algo:MatrixIRLS}
\begin{algorithmic}
\STATE{\bfseries Input:} Indices $\Omega$, observations $\f{y} \in \R^m$, rank estimate $\widetilde{r}$.
\STATE Initialize $k=0$, $\epsilon^{(0)}=\infty$ and $W^{(0)} = \Id$.

\FOR{$k=1$ to $K$}
\STATE \textbf{Solve weighted least squares:} Use a \emph{conjugate gradient method} to solve
\vspace*{-2mm}
\begin{equation} \label{eq:MatrixIRLS:Xdef}
\f{X}^{(k)} =\argmin\limits_{P_{\Omega}(\f{X})=\f{y}} \langle \f{X}, W^{(k-1)}(\f{X}) \rangle.
\end{equation}
\vspace*{-4mm}
\STATE \textbf{Update smoothing:} \label{eq:MatrixIRLS:bestapprox} Compute 	$\widetilde{r}+1$-th singular value of $\f{X}\hk$ to update
\vspace*{-5mm}
\begin{equation} \label{eq:MatrixIRLS:epsdef}
\epsilon_k=\min\left(\epsilon_{k-1},\sigma_{\widetilde{r}+1}(\f{X}\hk)\right).
\end{equation}
\vspace*{-3mm}

\STATE \textbf{Update weight operator:} For $r_k := |\{i \in [d]: \sigma_i(\f{X}\hk) > \epsilon_k\}|$, compute the first $r_k$ singular values $\sigma_i\hk := \sigma_i(\f{X}\hk)$ and matrices $\f{U}\hk \in \R^{d_1 \times r_k}$ and $\f{V}\hk \in \R^{d_2 \times r_k}$ with leading $r_k$ left/ right singular vectors of $\f{X}\hk$ to update $W\hk$ defined in \Cref{eq:def:W}. \label{eqdef:Wk} 

 \vspace{-.3cm}
\ENDFOR
\STATE{\bfseries Output:} $\f{X}^{(K)}$.
\end{algorithmic}
\end{algorithm}

\vspace{-4mm}
\section{Computational Complexity}
\vspace{-2mm}
A crucial property of \Cref{algo:MatrixIRLS} is that due to the choice of the smoothing function \cref{eq:smoothing:Fpeps}, the weight operator \cref{eq:def:W} and the smoothing update rule \cref{eq:MatrixIRLS:epsdef}, the action $W\hk$ on a matrix $\f{Z}$ can be implemented by a scalar multiplications and multiplications with rectangular $(d_1 \times r_k)$- and $(r_k \times d_2)$-matrices.
Utilizing this specific structure of the weight operators, we obtain an implementation of \texttt{MatrixIRLS} with a time and space complexity of the same order as for state-of-the-art first-order algorithms based on matrix factorization \cite{chen_chi18}. We refer to the supplementary materials for details and a proof.
\begin{theorem} \label{thm:MatrixIRLS:computationalcost:Xkk}
	Let $\f{X}\hk \in \Rdd$ be the $k$-th iterate of \texttt{MatrixIRLS} for an observation vector $\f{y} \in \R^m$ and $\widetilde{r}=r$. Assume that $\sigma_i\hk \leq \epsilon_k$ for all $i > r$ and $\sigma_r\hk > \epsilon_k$. Then an implicit representation of the new iterate $\f{X}\hkk \in \R^{d_1 \times d_2}$ can be calculated in a \emph{time complexity} of 
	\[
	 O \left( (m r + r^2 D) \cdot N_{\text{CG\_inner}} \right),
	\]
where $N_{\text{CG\_inner}}$ is the number of inner iterations used in the conjugate gradient method and $D=\max(d_1,d_2)$. More precisely, $\f{X}\hkk$ can be represented as
\[
\f{X}\hkk = P_{\Omega}^*(\f{r}_{k+1})  +\f{U}\hk \f{M}_{1}^{(k+1)*} + \f{M}_2^{(k+1)} \f{V}^{(k)*},
\]
where $\f{r}_{k+1} \in \R^m$, $\f{M}_{1}^{(k+1)} \in \R^{d_2 \times r}$ and $\f{M}_2^{(k+1)} \in \R^{d_1 \times r}$,
i.e., with a \emph{space complexity} of $O ( m+ r D)$.
\end{theorem}
\vspace{-2mm}
\Cref{thm:MatrixIRLS:computationalcost:Xkk} illustrates the computational advantage of \texttt{MatrixIRLS} compared to previous iteratively reweighted least squares algorithms for low-rank matrix recovery problems \cite{Fornasier11,Mohan10,KS18}, which all require the storage and updates of full $(d_1 \times d_2)$-matrices and the calculation of singular value decompositions of these. We comment on why a constant number of inner iterations $N_{\text{CG\_inner}}$ typically suffices in the supplementary materials.

As $P_{\Omega}^*(\f{r}_{k+1}) \in \Rdd$ is $m$-sparse, $\f{X}\hkk$ is a sum of a sparse and two rank-$r$ matrices and thus fast matrix-vector multiplications can be used in methods such as Lanczos bidiagonalization or randomized Block Krylov \cite{MuscoMusco15} to compute $r_{k+1}$ singular values and vectors of $\f{X}\hkk$ in step 3 of \Cref{algo:MatrixIRLS}. 
 \vspace*{-4mm}
\section{Theoretical Analysis} \label{sec:convergence}
\vspace{-2mm}
As it is beyond the scope of this format, we leave a detailed convergence analysis of \texttt{MatrixIRLS} to an upcoming paper. By establishing a global majorization property of the quadratic model function implicitly defined by the weight operator $W\hk$, it is possible to show that accumulation points of $(\f{X}\hk)_{k \geq 1}$ are stationary points of the $\overline{\epsilon}$-smoothed log-det objective $F_{\overline{\epsilon}}$ for $\overline{\epsilon}:=\lim_{k\to \infty} \epsilon_k$. Furthermore, we can establish local convergence \emph{with quadratic rate} of $\texttt{MatrixIRLS}$ to an incoherent rank-$r$ ground truth $\f{X}^0 \in \Rdd$ if a random sampling set $\Omega$ of size $\Omega(r(d_1 + d_2) \log(d_1 + d_2))$ is given, with high probability. In this result, the bound on the sample complexity does \emph{not} depend on the condition number of $\f{X}^0$. Such a result is new for matrix completion, as the results of \cite{KS18} studying a similar algorithm only cover measurement operators fulfilling a \emph{null space property}. 

\vspace*{-2.5mm}
\subsection{MatrixIRLS as saddle-escaping smoothing Newton method} \label{sec:Newton:interpretation}
From a theoretical point of view, the local quadratic convergence rate is an inherently local property that does not explain the numerically observed global convergence, which is remarkable due to the non-convexity of the objective.

A possible avenue to explain this is to interpret \texttt{MatrixIRLS} as a \emph{saddle-escaping smoothing Newton method}. Smoothing Newton methods minimize a non-smooth and possibly non-convex function $F$ by using derivatives of certain smoothings of $F$ \cite{ChenQiSun98,Chen2012smoothing}. Interpreting the optimization problem
$\min_{\f{X}:P_{\Omega}(X)=\f{y}} F_{\epsilon_k}(\f{X})$
as an unconstrained optimization problem over the null space of $P_{\Omega}$, we can write
\vspace{-.14cm}
\begin{equation*}
    \begin{split}
&\f{X}^{(k+1)} = \f{X}^{(k)} - P_{\Omega^c}^* \left(P_{\Omega^c} W\hk P_{\Omega^c}^*\right)^{-1} P_{\Omega^c} W\hk (\f{X}\hk) \\
\end{split}
\end{equation*}
\begin{equation*}
    \begin{split}
	&\!\!= \f{X}^{(k)}\!\! -\!P_{\Omega^c}^* \left(P_{\Omega^c}\!  \overline{\nabla^2 F_{\epsilon_k}(\f{X}\hk)}\! P_{\Omega^c}^*\right)^{-1}\!\!\!\!\!\! P_{\Omega^c} \nabla F_{\epsilon}(\f{X}\hk),
\end{split}
\end{equation*}
if $\Omega^c = [d_1] \times [d_2] \setminus \Omega$ corresponds to the unobserved indices, where $\overline{\nabla^2 F_{\epsilon_k}(\f{X}\hk)}: \Rdd \to \Rdd$ is a \emph{modified} Hessian of $F_{\epsilon_k}$ at $\f{X}\hk$ that replaces negative eigenvalues of the Hessian $\nabla^2 F_{\epsilon_k}(\f{X}\hk)$ by positive ones and slightly increases small eigenvalues. We refer to the supplementary material for more details. In \cite{PaternainMokhtariRibeiro19}, it has been proved that for a fixed smooth function $F_{\epsilon_k}$, similar modified Newton-type steps are able to escape the first-order saddle points at a rate that is independent of the problem's condition number.
\vspace*{-3.5mm}
\section{Numerical Experiments} \label{sec:numerics}
We explore the performance of \texttt{MatrixIRLS} for the completion of synthetic low-rank matrices in terms of data efficiency and computational efficiency in comparison to state-of-the-art algorithms in the literature. Algorithmic parameters are detailed in the supplementary materials. We consider the following setup: We sample a pair of random matrices $\f{U} \in \R^{d_1 \times r}$ and $\f{V} \in \R^{d_2 \times r}$ with $r$ orthonormal columns, and define the diagonal matrix $\Sigma \in \R^{r \times r}$ such that $\Sigma_{ii} = \kappa \exp(-\log(\kappa)\frac{i-1}{r-1})$ for $i \in [r]$. With this definition, we define a ground truth matrix $\f{X}^0 =\f{U}\Sigma \f{V}^*$ of rank $r$ that has exponentially decaying singular values between $\kappa$ and $1$. Furthermore, for a given factor $\rho \geq 1$, we sample a set $\Omega \subset [d_1] \times [d_2]$ of size $m = \lfloor \rho r (d_1 +d_2 - r) \rfloor$ indices randomly without replacement, implying that $\rho$ can be interpreted as an \emph{oversampling factor} as $r (d_1 +d_2 - r)$ is equal to the number of degrees of freedom of an $(d_1 \times d_2)$-dimensional rank-$r$ matrix.\footnote{The experiments can be reproduced with the code provided in \url{https://github.com/ckuemmerle/MatrixIRLS}.}

\vspace*{-2.5mm} 
\subsection{Data-efficient recovery of ill-conditioned matrices}
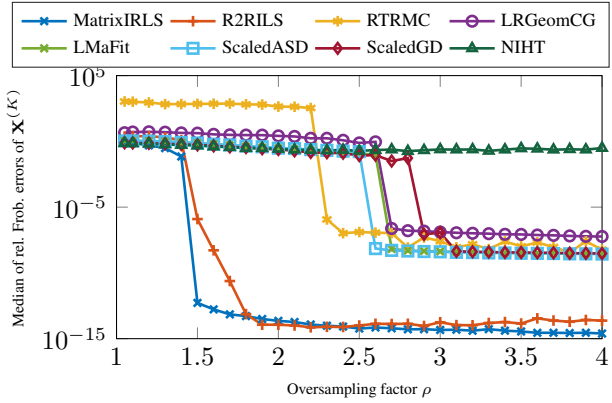
\begin{figure}[h!]
    \setlength\figureheight{35mm} 
    \setlength\figurewidth{66mm}
%
%
\definecolor{mycolor1}{rgb}{0.00000,0.44700,0.74100}
\definecolor{mycolor2}{rgb}{0.85000,0.32500,0.09800}
\definecolor{mycolor3}{rgb}{0.92900,0.69400,0.12500}
\definecolor{mycolor4}{rgb}{0.49400,0.18400,0.55600}
\definecolor{mycolor5}{rgb}{0.46600,0.67400,0.18800}
\definecolor{mycolor6}{rgb}{0.30100,0.74500,0.93300}
\definecolor{mycolor7}{rgb}{0.63500,0.07800,0.18400}
\definecolor{mycolor8}{rgb}{0.08000,0.39200,0.25100}

\begin{tikzpicture}

\begin{axis}[%
width=0.978\figurewidth,
height=\figureheight,
at={(0\figurewidth,0\figureheight)},
scale only axis,
xmin=1,
xmax=4,
xlabel style={font=\color{white!15!black}},
xlabel={Oversampling factor $\rho$},
ymode=log,
ymin=1e-15,
ymax=100000,
yminorticks=true,
ylabel style={font=\color{white!15!black}},
ylabel={Median of rel. Frob. errors of $\mathbf{X}^{(K)}$},
axis background/.style={fill=white},
legend style={legend cell align=left, align=left, draw=white!15!black},
xlabel style={font=\tiny},ylabel style={font=\tiny},legend style={font=\fontsize{7}{30}\selectfont, anchor=south, legend columns = 4, at={(0.4,1.03)}}
]
\addplot [color=mycolor1, line width=1.0pt, mark=x, mark options={solid, mycolor1}]
  table[row sep=crcr]{%
1.05	0.854255555628291\\
1.1	0.749423097099185\\
1.2	0.494716881980835\\
1.3	0.307031095563721\\
1.4	0.0710950413951957\\
1.5	5.22910462862914e-13\\
1.6	1.68588918367454e-13\\
1.7	7.0339665885088e-14\\
1.8	5.22646777788465e-14\\
1.9	3.12007321555135e-14\\
2	2.26159575881266e-14\\
2.1	1.84275579519234e-14\\
2.2	1.20010224139451e-14\\
2.3	9.65828095082859e-15\\
2.4	8.33489396563835e-15\\
2.5	5.87405834953726e-15\\
2.6	7.06122879837332e-15\\
2.7	6.140965219143e-15\\
2.8	5.15640770569221e-15\\
2.9	5.24064290362918e-15\\
3	4.58553528727644e-15\\
3.1	4.81251473302107e-15\\
3.2	4.03701051254568e-15\\
3.3	5.28365319386791e-15\\
3.4	4.07526625097761e-15\\
3.5	3.54126224418117e-15\\
3.6	2.79718832356394e-15\\
3.7	2.81718954270358e-15\\
3.8	2.69583211552255e-15\\
3.9	2.82440163639793e-15\\
4	2.46058369466678e-15\\
};
\addlegendentry{MatrixIRLS}

\addplot [color=mycolor2, line width=1.0pt, mark=+, mark options={solid, mycolor2}]
  table[row sep=crcr]{%
1.05	2.26456245685704\\
1.1	2.5024356292252\\
1.2	2.43235945106256\\
1.3	2.03349200850792\\
1.4	1.21611319874659\\
1.5	1.23751335798688e-06\\
1.6	5.20820599960089e-09\\
1.7	2.4080827102596e-11\\
1.8	7.35611767992492e-14\\
1.9	1.14972497705595e-14\\
2	1.14753038353458e-14\\
2.1	1.00726730331893e-14\\
2.2	7.34438613648422e-15\\
2.3	7.96384643655968e-15\\
2.4	7.79328724569884e-15\\
2.5	1.00231206491362e-14\\
2.6	1.3546830583003e-14\\
2.7	1.31005965452693e-14\\
2.8	1.35392160987976e-14\\
2.9	9.22728400705331e-15\\
3	1.77220221973004e-14\\
3.1	1.07494351508938e-14\\
3.2	1.04816978255403e-14\\
3.3	1.55757331921461e-14\\
3.4	1.83155751788827e-14\\
3.5	1.41801812290188e-14\\
3.6	3.49996967411273e-14\\
3.7	2.38941709817596e-14\\
3.8	1.93508288410887e-14\\
3.9	2.6480581760627e-14\\
4	2.35937791312582e-14\\
};
\addlegendentry{R2RILS}

\addplot [color=mycolor3, line width=1.0pt, mark=asterisk, mark options={solid, mycolor3}]
  table[row sep=crcr]{%
1.05	1075.26057752175\\
1.1	1036.88076801983\\
1.2	898.187957850437\\
1.3	676.778169415698\\
1.4	684.417922044455\\
1.5	711.785063703726\\
1.6	725.442000383892\\
1.7	751.035923142478\\
1.8	651.905082795613\\
1.9	629.404587607984\\
2	441.777176221888\\
2.1	429.821955824141\\
2.2	343.56802178822\\
2.3	1.15941912323249e-06\\
2.4	1.04551680418991e-07\\
2.5	1.27275200027277e-07\\
2.6	1.14870558919892e-07\\
2.7	1.05903565530647e-07\\
2.8	8.54610944390847e-09\\
2.9	4.96007271424229e-08\\
3	2.61470122950561e-08\\
3.1	8.86485982905778e-09\\
3.2	1.49412147034074e-08\\
3.3	6.46692958302814e-09\\
3.4	2.4700422542626e-08\\
3.5	1.03709991798648e-08\\
3.6	2.06415898459728e-08\\
3.7	9.18382240490974e-09\\
3.8	3.58229325687377e-09\\
3.9	2.52518262515746e-08\\
4	5.83585108699413e-09\\
};
\addlegendentry{RTRMC}

\addplot [color=mycolor4, line width=1.0pt, mark=o, mark options={solid, mycolor4}]
  table[row sep=crcr]{%
1.05	4.54218788203825\\
1.1	5.12790257479795\\
1.2	5.33019386122026\\
1.3	4.98053994216232\\
1.4	4.39615077981243\\
1.5	4.22488135311128\\
1.6	3.4144198382766\\
1.7	3.12003690973836\\
1.8	3.00097875387249\\
1.9	2.983446174957\\
2	2.60611156994507\\
2.1	2.36681710851521\\
2.2	1.71475814466035\\
2.3	1.71775747587247\\
2.4	1.20879362205861\\
2.5	0.764388428225149\\
2.6	0.91644120517575\\
2.7	2.39041820122551e-07\\
2.8	1.62071295437065e-07\\
2.9	1.50039926066069e-07\\
3	1.27835472461951e-07\\
3.1	1.12400902325457e-07\\
3.2	1.03200233643693e-07\\
3.3	9.61158837664088e-08\\
3.4	8.59690579442177e-08\\
3.5	8.40231752692976e-08\\
3.6	7.59751041476831e-08\\
3.7	7.16773678345403e-08\\
3.8	6.4853253141164e-08\\
3.9	6.01399241139154e-08\\
4	5.84948579668478e-08\\
};
\addlegendentry{LRGeomCG}

\addplot [color=mycolor5, line width=1.0pt, mark=x, mark options={solid, mycolor5}]
  table[row sep=crcr]{%
1.05	1.12746309687931\\
1.1	1.1215115543341\\
1.2	1.05005289952185\\
1.3	0.979491896633853\\
1.4	0.827263202581372\\
1.5	0.738733099621636\\
1.6	0.607355931786455\\
1.7	0.497511329328573\\
1.8	0.405136304159793\\
1.9	0.334586407456963\\
2	0.283213525406684\\
2.1	0.242367443226484\\
2.2	0.196034213853877\\
2.3	0.151317020176626\\
2.4	0.117156714237261\\
2.5	0.115895091936267\\
2.6	0.101974980837135\\
2.7	6.65749401699423e-09\\
2.8	5.01786293007783e-09\\
2.9	4.59006342081414e-09\\
3	4.28262136377922e-09\\
3.1	3.97593246298878e-09\\
3.2	3.84862100638697e-09\\
3.3	3.56396882136211e-09\\
3.4	3.44257715382382e-09\\
3.5	3.28634578928152e-09\\
3.6	3.22147596465297e-09\\
3.7	3.1379469525384e-09\\
3.8	2.97955273363032e-09\\
3.9	2.94591484325744e-09\\
4	2.8687371350171e-09\\
};
\addlegendentry{LMaFit}

\addplot [color=mycolor6, line width=1.0pt, mark=square, mark options={solid, mycolor6}]
  table[row sep=crcr]{%
1.05	1.06466416015086\\
1.1	1.06463314388243\\
1.2	1.03965187979647\\
1.3	0.960668507770738\\
1.4	0.824747589865423\\
1.5	0.745062952365437\\
1.6	0.604212717779343\\
1.7	0.533104871209988\\
1.8	0.434557035092144\\
1.9	0.388395943543961\\
2	0.305635793463275\\
2.1	0.284434146621849\\
2.2	0.235606018890938\\
2.3	0.208483454008394\\
2.4	0.182529631516212\\
2.5	0.165489912583052\\
2.6	7.05991430201768e-09\\
2.7	5.16997427925381e-09\\
2.8	4.59969613371997e-09\\
2.9	4.21095337333196e-09\\
3	4.01020980931519e-09\\
3.1	3.65962492755945e-09\\
3.2	3.51540016368054e-09\\
3.3	3.30205204705014e-09\\
3.4	3.24585144342233e-09\\
3.5	3.12736215079214e-09\\
3.6	3.06860576382731e-09\\
3.7	2.98223020227907e-09\\
3.8	2.81961243454812e-09\\
3.9	2.81855322114525e-09\\
4	2.70286224887889e-09\\
};
\addlegendentry{ScaledASD}

\addplot [color=mycolor7, line width=1.0pt, mark=diamond, mark options={solid, mycolor7}]
  table[row sep=crcr]{%
1.05	0.780342983376787\\
1.1	0.760559629534551\\
1.2	0.716937289667204\\
1.3	0.655527483255922\\
1.4	0.589124333833299\\
1.5	0.542132854431097\\
1.6	0.45418458358205\\
1.7	0.385560512899637\\
1.8	0.341913628226789\\
1.9	0.278746040604984\\
2	0.242730595784693\\
2.1	0.192126001507383\\
2.2	0.145452661065458\\
2.3	0.138749643187959\\
2.4	0.127851513423253\\
2.5	0.0834717206401568\\
2.6	0.0896687668276914\\
2.7	0.0318824695108248\\
2.8	0.0534563334227891\\
2.9	8.57397358047281e-08\\
3	1.1479597956224e-07\\
3.1	4.3746445456053e-09\\
3.2	4.05807191058956e-09\\
3.3	3.86537463377338e-09\\
3.4	3.63303397873685e-09\\
3.5	3.57399045753236e-09\\
3.6	3.43987811253563e-09\\
3.7	3.34572184308354e-09\\
3.8	3.1544837589476e-09\\
3.9	3.12171693602334e-09\\
4	2.97912161408277e-09\\
};
\addlegendentry{ScaledGD}

\addplot [color=mycolor8, line width=1.0pt, mark=triangle, mark options={solid, mycolor8}]
  table[row sep=crcr]{%
1.05	0.741178657713085\\
1.1	0.710345115228534\\
1.2	0.657484763497934\\
1.3	0.583841951688662\\
1.4	0.52969749416985\\
1.5	0.476354734852938\\
1.6	0.404884442048636\\
1.7	0.360820981001756\\
1.8	0.306935251980735\\
1.9	0.303708974319308\\
2	0.267456591481581\\
2.1	0.234331899499627\\
2.2	0.21993357166271\\
2.3	0.206373745234411\\
2.4	0.224257420213961\\
2.5	0.188961454098439\\
2.6	0.21827413665062\\
2.7	0.238458775410865\\
2.8	0.180849453808525\\
2.9	0.211087228946457\\
3	0.245008088121678\\
3.1	0.236380488132301\\
3.2	0.243076273505301\\
3.3	0.189558676300446\\
3.4	0.234927232199093\\
3.5	0.291938223692077\\
3.6	0.274422676296374\\
3.7	0.233688337092644\\
3.8	0.250113008520314\\
3.9	0.235759425673383\\
4	0.312027766249649\\
};
\addlegendentry{NIHT}

\end{axis}
\end{tikzpicture}%
\vspace*{-8mm}
\caption{Comparison of matrix completion algorithms for $1000 \times 1000$ matrices of rank $r=5$ with condition number $\kappa=10$, given $m= \lfloor \rho r (d_1 + d_2-r)\rfloor$ random samples. Median of Frobenius errors $\|\f{X}^{(K)}-\f{X}^0\|_F/\|\f{X}^0\|_F$ of $100$ independent realizations.}
\label{fig:sampcomp:1}
\end{figure}

First, we run \texttt{MatrixIRLS} and the algorithms \texttt{R2RILS} \cite{BauchNadler20}, \texttt{RTRMC} \cite{boumal_absil_15}, \texttt{LRGeomCG} \cite{Vandereycken13}, \texttt{LMaFit} \cite{Wen12}, \texttt{ScaledASD} \cite{TannerWei16}, \texttt{ScaledGD} \cite{tong_ma_chi} and \texttt{NIHT} \cite{TannerWei13} to complete $\f{X}^0$ from $P_{\Omega}(\f{X}^0)$ where $\Omega$ corresponds to different oversampling factors $\rho$ between $1$ and $4$, and where the condition number of $\f{X}^0$ is $\kappa= \sigma_1(\f{X}^0)/\sigma_r(\f{X}^0) = 10$. In \Cref{fig:sampcomp:1}, we report the median Frobenius errors $\|\f{X}^{(K)}-\f{X}^0\|_F/\|\f{X}^0\|_F$ of the the respective algorithmic outputs $\f{X}^{(K)}$ across $100 $ independent realizations. 

\vspace{-.3cm}
We see that \texttt{MatrixIRLS} and \texttt{R2RILS} are the only algorithms that are able to complete $\f{X}^0$ already for $\rho=1.5$, whereas the other algorithms, except from \texttt{NIHT}, are able to reconstruct the matrix most of the times if $\rho$ is at least between $2.4$ and $3.0$. This confirms the findings of \cite{BauchNadler20} that even for quite well-conditioned matrices, fewer samples are required if second-order methods such as \texttt{R2RILS} or \texttt{MatrixIRLS} are used.

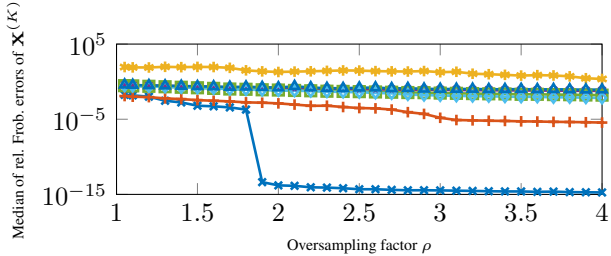
\begin{figure}[ht]
    \setlength\figureheight{25mm} 
    \setlength\figurewidth{66mm}
%
%

\definecolor{mycolor1}{rgb}{0.00000,0.44700,0.74100}
\definecolor{mycolor2}{rgb}{0.85000,0.32500,0.09800}
\definecolor{mycolor3}{rgb}{0.92900,0.69400,0.12500}
\definecolor{mycolor4}{rgb}{0.49400,0.18400,0.55600}
\definecolor{mycolor5}{rgb}{0.46600,0.67400,0.18800}
\definecolor{mycolor6}{rgb}{0.30100,0.74500,0.93300}
\definecolor{mycolor7}{rgb}{0.63500,0.07800,0.18400}
\definecolor{mycolor8}{rgb}{0.08000,0.39200,0.25100}
\begin{tikzpicture}

\begin{axis}[%
width=0.978\figurewidth,
height=0.8\figureheight,
at={(0\figurewidth,0\figureheight)},
scale only axis,
xmin=1,
xmax=4,
xlabel style={font=\color{white!15!black}},
xlabel={Oversampling factor $\rho$},
ymode=log,
ymin=1e-15,
ymax=100000,
yminorticks=true,
ylabel style={font=\color{white!15!black}},
ylabel={Median of rel. Frob. errors of $\mathbf{X}^{(K)}$},
axis background/.style={fill=white},
xlabel style={font=\tiny},ylabel style={font=\tiny},
]
\addplot [color=mycolor1, line width=1.0pt, mark=x, mark options={solid, mycolor1}, forget plot]
  table[row sep=crcr]{%
1.05	0.017696793293338\\
1.1	0.0149247467929977\\
1.2	0.0100358741421393\\
1.3	0.00287401737080599\\
1.4	0.0017217492213112\\
1.5	0.000627509308504794\\
1.6	0.000523487605715263\\
1.7	0.000374934439857983\\
1.8	0.000192177543449912\\
1.9	4.26148893627606e-14\\
2	1.58623861308894e-14\\
2.1	1.37788927229948e-14\\
2.2	8.82670925448159e-15\\
2.3	8.02507277531667e-15\\
2.4	6.59203757233433e-15\\
2.5	4.85856656140693e-15\\
2.6	4.84646205110599e-15\\
2.7	4.03688126061311e-15\\
2.8	3.52166740720263e-15\\
2.9	3.57778726155854e-15\\
3	3.39853128746517e-15\\
3.1	3.13012670512855e-15\\
3.2	2.90704826967485e-15\\
3.3	2.64447657593206e-15\\
3.4	2.47776651290129e-15\\
3.5	2.47702997773112e-15\\
3.6	2.18614791909683e-15\\
3.7	2.35729374885045e-15\\
3.8	2.0121153940015e-15\\
3.9	1.92108224894144e-15\\
4	1.82116757856728e-15\\
};
\addplot [color=mycolor2, line width=1.0pt, mark=+, mark options={solid, mycolor2}, forget plot]
  table[row sep=crcr]{%
1.05	0.0121729439043043\\
1.1	0.00997144484447369\\
1.2	0.00872449609984643\\
1.3	0.00594736215558571\\
1.4	0.00476375327613421\\
1.5	0.00350287893314341\\
1.6	0.00286693482596602\\
1.7	0.00235681936745234\\
1.8	0.00176184450881462\\
1.9	0.00157710691467946\\
2	0.0012894886519888\\
2.1	0.000890416387229741\\
2.2	0.000612510871767183\\
2.3	0.000628320352780005\\
2.4	0.000422461566296125\\
2.5	0.000307890100903442\\
2.6	0.000278390558878799\\
2.7	0.000193359962186809\\
2.8	8.73876008786447e-05\\
2.9	5.12954109819375e-05\\
3	1.55842407560365e-05\\
3.1	7.8624565679366e-06\\
3.2	7.09562671875953e-06\\
3.3	6.33123971806548e-06\\
3.4	5.92027738932675e-06\\
3.5	5.24569952214572e-06\\
3.6	4.96587249559889e-06\\
3.7	4.51427820997321e-06\\
3.8	4.28391740187636e-06\\
3.9	3.97622888670062e-06\\
4	3.61094465249887e-06\\
};
\addplot [color=mycolor3, line width=1.0pt, mark=asterisk, mark options={solid, mycolor3}, forget plot]
  table[row sep=crcr]{%
1.05	88.1773642405158\\
1.1	79.3512035402019\\
1.2	76.7745266497778\\
1.3	94.3732403850792\\
1.4	84.7556749644728\\
1.5	92.0131576237411\\
1.6	88.5243100930591\\
1.7	68.9512284446234\\
1.8	32.1455731454929\\
1.9	23.4058432834348\\
2	19.4557226288538\\
2.1	21.6674555929097\\
2.2	23.8182185283779\\
2.3	27.0991116333954\\
2.4	31.168987747205\\
2.5	29.8500776136293\\
2.6	25.0704730678691\\
2.7	23.4754271963616\\
2.8	22.1912782891398\\
2.9	25.7075375758748\\
3	18.1710204054141\\
3.1	20.7392540885847\\
3.2	11.3573930075962\\
3.3	9.77324309823878\\
3.4	8.03811839649874\\
3.5	6.38279842759175\\
3.6	7.36222030845832\\
3.7	6.69451532674315\\
3.8	5.09370937566492\\
3.9	2.83881370325447\\
4	2.28326868894603\\
};
\addplot [color=mycolor4, line width=1.0pt, mark=x, mark options={solid, mycolor4}, forget plot]
  table[row sep=crcr]{%
1.05	0.370479094334638\\
1.1	0.344738688789208\\
1.2	0.305003685693449\\
1.3	0.274819969474776\\
1.4	0.247740013300637\\
1.5	0.220229797221742\\
1.6	0.197593429507371\\
1.7	0.174974560110361\\
1.8	0.161880928811185\\
1.9	0.149072956986771\\
2	0.135573313964275\\
2.1	0.125617559765276\\
2.2	0.117624754306401\\
2.3	0.107770796023168\\
2.4	0.102869148200105\\
2.5	0.0973699783171384\\
2.6	0.0911976587510287\\
2.7	0.087070467542029\\
2.8	0.0851695140911353\\
2.9	0.0820343098805306\\
3	0.0784133307659842\\
3.1	0.0768336670990977\\
3.2	0.073298345272641\\
3.3	0.071935717650967\\
3.4	0.0688059674895942\\
3.5	0.0675089237185283\\
3.6	0.0660170864969936\\
3.7	0.0630233750193977\\
3.8	0.0608667400397613\\
3.9	0.0599704565136605\\
4	0.0594213595172383\\
};
\addplot [color=mycolor5, line width=1.0pt, mark=square, mark options={solid, mycolor5}, forget plot]
  table[row sep=crcr]{%
1.05	0.293535194670914\\
1.1	0.273424548894986\\
1.2	0.241863789436971\\
1.3	0.210308160885623\\
1.4	0.179671150516937\\
1.5	0.153704204975937\\
1.6	0.127383037107324\\
1.7	0.104581966212891\\
1.8	0.0967381965419033\\
1.9	0.0829217878720201\\
2	0.0767587856465289\\
2.1	0.0677358836206883\\
2.2	0.0610941507309169\\
2.3	0.0536826770846946\\
2.4	0.0501433886034221\\
2.5	0.047032656931879\\
2.6	0.0424026045167179\\
2.7	0.0391535053499478\\
2.8	0.0355206008825965\\
2.9	0.0342638696279479\\
3	0.0301352181009652\\
3.1	0.0279617109248466\\
3.2	0.0257797090099662\\
3.3	0.0240692916355119\\
3.4	0.023210471331287\\
3.5	0.0206607719407823\\
3.6	0.0208466572877179\\
3.7	0.0190403180574802\\
3.8	0.0184171481393597\\
3.9	0.0167656690800692\\
4	0.0160719495980064\\
};
\addplot [color=mycolor6, line width=1.0pt, mark=diamond, mark options={solid, mycolor6}, forget plot]
  table[row sep=crcr]{%
1.05	0.282005628998929\\
1.1	0.259057103039847\\
1.2	0.222535496634528\\
1.3	0.200741125186381\\
1.4	0.17154972703974\\
1.5	0.151177149976177\\
1.6	0.127593953591546\\
1.7	0.103040137529032\\
1.8	0.0888678232426331\\
1.9	0.0775871141116796\\
2	0.0689565811530352\\
2.1	0.0628166001450546\\
2.2	0.0552483658972788\\
2.3	0.0479359822321287\\
2.4	0.0423040634728772\\
2.5	0.0392512972817216\\
2.6	0.0318814987640092\\
2.7	0.0274905957643751\\
2.8	0.0195835655807152\\
2.9	0.0161063827691886\\
3	0.0145834490638091\\
3.1	0.0128517676575462\\
3.2	0.0118565862006724\\
3.3	0.0113175533060402\\
3.4	0.0101385628631579\\
3.5	0.00956005369402189\\
3.6	0.0095480484530076\\
3.7	0.00816374445443083\\
3.8	0.00815332450461049\\
3.9	0.00738917286243342\\
4	0.00685805936868511\\
};
\addplot [color=mycolor1, line width=1.0pt, mark=triangle, mark options={solid, mycolor1}, forget plot]
  table[row sep=crcr]{%
1.05	0.381684333532095\\
1.1	0.350455776498883\\
1.2	0.31465709602788\\
1.3	0.289284561868943\\
1.4	0.261387299429177\\
1.5	0.236763350860108\\
1.6	0.219918434274913\\
1.7	0.193810606899661\\
1.8	0.209297075614819\\
1.9	0.181745589813183\\
2	0.163995242516317\\
2.1	0.24155823327221\\
2.2	0.164789720459069\\
2.3	0.136348972794721\\
2.4	0.131894690711273\\
2.5	0.154017219048052\\
2.6	0.159725873265842\\
2.7	0.127157617918866\\
2.8	0.132809280412949\\
2.9	0.110726691809945\\
3	0.103675178074836\\
3.1	0.105090353370288\\
3.2	0.105795441660664\\
3.3	0.0954892884532374\\
3.4	0.0971222604079868\\
3.5	0.100048891056211\\
3.6	0.0925294362854234\\
3.7	0.0961592579970626\\
3.8	0.0913540561972631\\
3.9	0.0848092286768069\\
4	0.0876188156232948\\
};
\end{axis}
\end{tikzpicture}%
\vspace*{-6mm}
\caption{Comparison of matrix completion algorithms as in \Cref{fig:sampcomp:1}, but with $\kappa = 10^5$. Median of $50$ realizations.}
\label{fig:sampcomp:2}
\end{figure}

We repeat this experiment for ill-conditioned matrices $\f{X}^0$ with $\kappa = 10^5$. In  \Cref{fig:sampcomp:2}, we see that current state-of-the-art methods are \emph{not able} to achieve exact recovery of $\f{X}^0$. This is in particular true as given the exponential decay of the singular values, in order to recover the subspace corresponding to the smallest singular value of $\f{X}^0$, a relative Frobenius error of $10^{-5}$ or even several orders of magnitude smaller needs to be achieved. We observe that \texttt{MatrixIRLS} is the only method that is able to complete $\f{X}^0$ for any of the considered oversampling factors.

\vspace*{-0.5cm}
 
\vspace{-.5cm}

\subsection{Running time for ill-conditioned problems}

\begin{figure}[ht]
    \setlength\figureheight{25mm} 
    \setlength\figurewidth{66mm}
\input{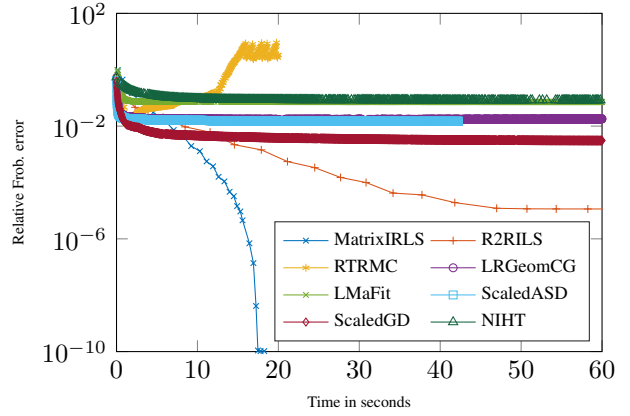}
\vspace*{-5mm}
\caption{Completion task for a highly ill-conditioned $1000 \times 1000$ matrix of rank $r=10$ with $\kappa=10^5$ ($\rho= 4$).}
\label{running_time_ill_cond}
\end{figure}

In Figure \ref{running_time_ill_cond}, for an oversampling rate of $\rho = 4$, we illustrate the completion of one single highly ill-conditioned matrix. We again can see that only second-order methods such as \texttt{R2RILS} or \texttt{MatrixIRLS} are able to achieve a relative Frobenius error $\approx 10^{-5}$ or smaller but \texttt{MatrixIRLS} does it in a much \emph{faster} way. Besides that, it also retrieves all the singular values with high precision by attaining a very low relative error.

\begin{figure}[h!]
    \setlength\figureheight{25mm} 
    \setlength\figurewidth{66mm}
%
%
\definecolor{mycolor1}{rgb}{0.00000,0.44700,0.74100}
\definecolor{mycolor2}{rgb}{0.85000,0.32500,0.09800}
\definecolor{mycolor3}{rgb}{0.92900,0.69400,0.12500}
\definecolor{mycolor4}{rgb}{0.49400,0.18400,0.55600}
\definecolor{mycolor5}{rgb}{0.46600,0.67400,0.18800}
\definecolor{mycolor6}{rgb}{0.30100,0.74500,0.93300}
\definecolor{mycolor7}{rgb}{0.63500,0.07800,0.18400}
\definecolor{mycolor8}{rgb}{0.08000,0.39200,0.25100}
\begin{tikzpicture}

\begin{axis}[%
width=0.978\figurewidth,
height=\figureheight,
at={(0\figurewidth,0\figureheight)},
scale only axis,
xmin=100,
xmax=1000,
xlabel={m (with n=m+100)},
ymin=2.5,
ymax=7,
ylabel={Execution time (seconds)},
axis background/.style={fill=white},
xmajorgrids,
ymajorgrids,
xlabel style={font=\tiny},ylabel style={font=\tiny},legend style={font=\fontsize{7}{30}\selectfont, anchor=south, legend columns = 2, at={(0.5,1.01)}}
]

(2,3) node {\tiny MatrixIRLS $r=5$}
\node [anchor=center] (MatrixIRLS) at (-20,-3) {\tiny MatrixIRLS $r=5$};
\addplot [color=mycolor1, line width=1pt,  mark=asterisk, mark options={solid, mycolor1}]
  table[row sep=crcr]{%
100	2.95304825418182\\
200	2.98230332166667\\
300	3.00724506338462\\
400	3.08807303557143\\
500	3.22394640293334\\
600	3.274575889875\\
700	3.33046376941177\\
800	3.42672688188889\\
900	3.54640258852632\\
1000	3.7073854512\\
};
\addlegendentry{MatrixIRLS r=5}

\addplot [color=mycolor2, line width=1pt,  mark=asterisk, mark options={solid, mycolor2}]
  table[row sep=crcr]{%
100	5.199734626\\
200	4.97052347\\
300	5.06955289830769\\
400	5.03849069399999\\
500	5.082422796\\
600	5.218918925\\
700	5.35722198082353\\
800	5.58250735188888\\
900	5.83916094263158\\
1000	6.4279443572\\
};
\addlegendentry{R2RILS r=5}

\addplot [color=mycolor1, line width=1.5pt]
  table[row sep=crcr]{%
100	3.55308134704762\\
200	3.44113479809091\\
300	3.37757395730435\\
400	3.35926170525\\
500	3.39726942392001\\
600	3.47673567853847\\
700	3.60797207814816\\
800	3.75406394928572\\
900	3.94470405496553\\
1000	4.20021874420001\\
};
\addlegendentry{MatrixIRLS r=10}

\addplot [color=mycolor2, line width=1.5pt]
  table[row sep=crcr]{%
100	6.20542460038095\\
200	6.087665208\\
300	6.03986127608695\\
400	6.08141432908333\\
500	6.20602667224\\
600	6.37619482830769\\
700	6.64026641185185\\
800	6.98286348607142\\
900	7.3768662233793\\
1000	7.76715193293332\\
};
\addlegendentry{R2RILS r=10}

\end{axis}
\end{tikzpicture}%
\vspace*{-5mm}
\caption{Comparison of \texttt{R2RILS} and \texttt{MatrixIRLS} for completion  of rank $r \in \{5,10\}$ matrices of size $m \times (m+100)$ and condition number $\kappa=10^2$ in terms of their execution time. Every point represents the average over 50 experiments. The other algorithms are not shown in this experiment because they typically do not reach a relative error below $10^{-4}$ for $\kappa\gg10^2$.}
\label{running_time_MatrixIRLS_vs_R2RILS}
\end{figure}
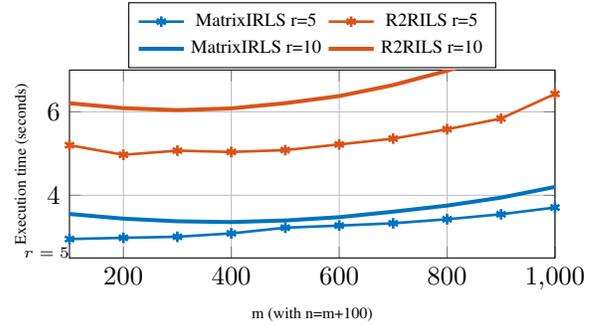

In Figure \ref{running_time_MatrixIRLS_vs_R2RILS}, we compare the execution time of \texttt{R2RILS} and \texttt{MatrixIRLS} when the dimension increases, for an oversampling rate of $\rho = 2.5$ with singular values linearly interpolated between $\kappa$ and $1$. We observe that the larger the dimensions are, the larger is the discrepancy in the running time of the two algorithms. This is accentuated for even higher condition numbers.

\vspace{-.8cm}

\section{Conclusion}\label{conclusion}
\vspace{-.2cm}
We formulated \texttt{MatrixIRLS}, a second order method that is able to efficiently complete large, highly ill-conditioned matrices from few samples, a problem for which most state-of-the-art methods fail. It improves on previous approaches for the optimization of non-convex rank objectives by applying a suitable smoothing strategy combined with saddle-escaping Newton-type steps.   

As one goal of our investigation has been also to provide an efficient implementation, we focused on the matrix completion problem, leaving the extension of the ideas to other low-rank matrix estimation problems to future work including the case of inexact data or measurement errors. Furthermore, while we establish a local convergence guarantee for the algorithm, a precise analysis of its global convergence behavior might be of interest.

\section*{Acknowledgements}
The authors thank Cong Ma for providing their code for the algorithm ScaledGD. The authors gratefully acknowledge the Leibniz Supercomputing Centre (LRZ) for providing computing time on their Linux cluster.  CMV gratefully acknowledges support by German Science Foundation (DFG) within the Gottfried Wilhelm Leibniz Prize under Grant BO 1734/20-1, under grant KR 4512/1-1, under contract number PO-1347/3-2 and within Germany's Excellence Strategy EXC-2111 390814868.

\bibliography{MatrixIRLS_wellconditioned}
\bibliographystyle{icml2020}

\newpage
\onecolumn

\appendix

\icmltitle{Supplementary material for \emph{Escaping Saddle Points in Ill-Conditioned Matrix Completion with a Scalable Second Order Method}}

This supplementary material is divided into three parts: In \Cref{sec:proof:computational}, we provide details on the implementation of \texttt{MatrixIRLS} and show \Cref{thm:MatrixIRLS:computationalcost:Xkk}. In \Cref{sec:remarks:smoothing:newton}, we provide some details about the interpretation of \texttt{MatrixIRLS} as a saddle-escaping smoothing Newton method, cf. \Cref{sec:Newton:interpretation}. Finally, we present a detailed description of the algorithmic parameters of the experiments of \Cref{sec:numerics} in \Cref{sec:experimental:details}, as well as one additional experiment.

\section{Proof of \Cref{thm:MatrixIRLS:computationalcost:Xkk}} \label{sec:proof:computational}
We first show how the action of the weight operator $W\hk: \Rdd \to \Rdd$ of \texttt{MatrixIRLS} can be represented just by scalar multiplcations and by multiplications by very rectangular matrices whose smaller dimension is $r_k := |\{i \in [d]: \sigma_i(\f{X}\hk) > \epsilon_k\}| = |\{i \in [d]: \sigma_i\hk > \epsilon_k\}|$.

To show this, let 
\begin{equation} \label{eq:Xk:bothsvds}
\f{X}\hk= \f{U}_k \dg(\sigma^{(k)}) \f{V}_k^* = 
\begin{bmatrix} 
\f{U}\hk & \f{U}_{\perp}\hk
\end{bmatrix}\begin{bmatrix} 
\f{\Sigma}\hk & 0 \\
0 & \f{\Sigma}_{\perp}\hk 
\end{bmatrix}
\begin{bmatrix} 
\f{V}^{(k)*} \\
\f{V}_{\perp}^{(k)*} 
\end{bmatrix}
\end{equation}
be a singular value decomposition of $\f{X}\hk$ with orthonormal matrices $\f{U}_k \in \R^{d_1 \times d_1}$ and $\f{V}_k \in \R^{d_2 \times d_2}$, and corresponding submatrices $\f{U}\hk \in \R^{d_1 \times r_k}$, $\f{U}_{\perp}\hk \in \R^{d_1 \times (d_1 -r_k)}$, $\f{V}\hk \in \R^{d_2 \times r_k}$, $\f{V}_{\perp}\hk \in \R^{d_2 \times (d_2 -r_k)}$, $\f{\Sigma}\hk := \diag(\sigma_1\hk,\ldots \sigma_{r_k}\hk)$ and $\f{\Sigma}_{\perp}\hk := \dg(\sigma_{r_k+1}\hk,\ldots \sigma_{d}\hk)$. Using this and the definition \cref{eq:def:W} of $W\hk$ and the matrix $\f{H}_k \in \Rdd$, it holds that for any $\f{Z} \in \Rdd$, 
\begin{equation*} 
\begin{split}
& W\hk(\f{Z}) = \f{U}_k \left[\f{H}_k \circ (\f{U}_k^{*} \f{Z} \f{V}_k)\right] \f{V}_k^{*} \\
&=\begin{bmatrix} 
	\f{U}\hk & \f{U}_{\perp}\hk
\end{bmatrix}
\left(\f{H}_k
\circ
\begin{bmatrix}
\f{U}^{(k)*} \f{Z} \f{V}^{(k)} &  \f{U}^{(k)*} \f{Z} \f{V}_{\perp}^{(k)} \\
\f{U}_{\perp}^{(k)*} \f{Z} \f{V}^{(k)} & \f{U}_{\perp}^{(k)*} \f{Z} \f{V}_{\perp}^{(k)} 
\end{bmatrix}
\right)
\begin{bmatrix} 
	\f{V}^{(k)*} \\ \f{V}_{\perp}^{(k)*}
\end{bmatrix} \\
&= \begin{bmatrix} 
	\f{U}\hk & \f{U}_{\perp}\hk
\end{bmatrix}
\left(
\begin{bmatrix}
\f{H}\hk &  \f{H}_{1,2}\hk \\
\f{H}_{2,1}\hk & \epsilon_k^{-2} \mathbf{1}
\end{bmatrix} 
\circ 
\begin{bmatrix}
\f{U}^{(k)*} \f{Z} \f{V}^{(k)} &  \f{U}^{(k)*} \f{Z} \f{V}_{\perp}^{(k)} \\
\f{U}_{\perp}^{(k)*} \f{Z} \f{V}^{(k)} & \f{U}_{\perp}^{(k)*} \f{Z} \f{V}_{\perp}^{(k)} 
\end{bmatrix}
\right)
\begin{bmatrix} 
	\f{V}^{(k)*} \\ \f{V}_{\perp}^{(k)*}
\end{bmatrix}  \\
&=
\begin{bmatrix} 
	\f{U}\hk & \f{U}_{\perp}\hk
\end{bmatrix}
\begin{bmatrix}
\f{H}\hk \circ \left(\f{U}^{(k)*} \f{Z} \f{V}^{(k)}\right) &  \f{D}\hk \f{U}^{(k)*} \f{Z} \f{V}_{\perp}^{(k)} \\
\f{U}_{\perp}^{(k)*} \f{Z} \f{V}^{(k)}\f{D}\hk  & \epsilon_k^{-2} \f{U}_{\perp}^{(k)*} \f{Z} \f{V}_{\perp}^{(k)} 
\end{bmatrix}
\begin{bmatrix} 
	\f{V}^{(k)*} \\ \f{V}_{\perp}^{(k)*}
\end{bmatrix},
\end{split}
\end{equation*}
dividing $\f{H}_k$ into the four block matrices 
\begin{equation} \label{eq:Hk:blockstructure}
\f{H}_k=
\left[
\begin{array}{c|c}
\f{H}\hk & \f{H}_{1,2}\hk \\
\hline
\f{H}_{2,1}\hk & \epsilon_k^{-2} \mathbf{1}
\end{array}
\right]
\end{equation}
such that $\f{H}\hk \in \R^{r_k \times r_k}$ with
\begin{equation} \label{eq:H:def:small}
\f{H}_{ij}\hk = \left(\sigma_i^{(k)} \sigma_j^{(k)}\right)^{-1}
\end{equation}
for all $i,j \in [r_k]$, $\f{H}_{1,2}\hk \in \R^{r_k \times (d_2 - r_k)}$ with $\left(\f{H}_{1,2}\hk\right)_{ij} =
\left(\sigma_i^{(k)} \epsilon_k\right)^{-1}$ for all $i \in [r_k]$ and $j \in [d_2-r_k]$, $\f{H}_{2,1}\hk \in \R^{ (d_1 - r_k) \times r_k}$ such that 
$ \left(\f{H}_{2,1}\hk\right)_{ij} = \left(\epsilon_k \sigma_j^{(k)} \right)^{-1}$ 
for all $i \in [d_1-r_k]$ and $j \in [r_k]$, and the $((d_1 - r_k) \times (d_2 -r_k))$-matrix of ones $\mathbf{1}$, and lastly, the diagonal matrix $\f{D}\hk \in \R^{r_k \times r_k}$ with
\begin{equation} \label{eq:D:def:small}
\f{D}_{ii}\hk = \left(\sigma_i^{(k)} \epsilon_k\right)^{-1}
\end{equation}
for all $i \in [r_k]$. Simplifying the last equality, we obtain then 
	\begin{equation} \label{eq:W:simplification}
	\begin{split}
W\hk(\f{Z}) &= \f{U}\hk \left[ \f{H}\hk \circ (\f{U}^{(k)*} \f{Z} \f{V}^{(k)})\right]\f{V}^{(k)*} + \f{U}\hk \f{D}\hk \f{U}^{(k)*}\f{Z}\f{V}_{\perp}^{(k)}\f{V}_{\perp}^{(k)*} \\
&+  \f{U}_{\perp}^{(k)}\f{U}_{\perp}^{(k)*} \f{Z} \f{V}\hk \f{D}\hk \f{V}^{(k)*} + \epsilon_k^{-2} \f{U}_{\perp}^{(k)}\f{U}_{\perp}^{(k)*} \f{Z}\f{V}_{\perp}^{(k)}\f{V}_{\perp}^{(k)*}  \\
&= \f{U}\hk 
	\left[ \f{H}\hk \circ (\f{U}^{(k)*} \f{Z} \f{V}^{(k)})\right]\f{V}^{(k)*} + \f{U}\hk \f{D}\hk \f{U}^{(k)*}\f{Z}(\f{I}-\f{V}^{(k)}\f{V}^{(k)*}) \\
&+  (\f{I}-\f{U}^{(k)}\f{U}^{(k)*}) \f{Z} \f{V}\hk \f{D}\hk \f{V}^{(k)*} + \epsilon_k^{-2} \left(\f{I}-\f{U}\hk \f{U}^{(k)*}\right) \f{Z} \left(\f{I}-\f{V}\hk \f{V}^{(k)*}\right) \\
&= \f{U}\hk 
	\left[ \left(\f{H}\hk- \epsilon_k^{-2}\mathbf{1}\right) \circ (\f{U}^{(k)*} \f{Z} \f{V}^{(k)})\right]\f{V}^{(k)*} + \f{U}\hk \left(\f{D}\hk- \epsilon_k^{-2} \f{I}\right) \f{U}^{(k)*}\f{Z}(\f{I}-\f{V}^{(k)}\f{V}^{(k)*}) \\
&+  (\f{I}-\f{U}^{(k)}\f{U}^{(k)*}) \f{Z} \f{V}\hk \left(\f{D}\hk- \epsilon_k^{-2} \f{I}\right) \f{V}^{(k)*} + \epsilon_k^{-2} \f{Z},
    \end{split}	
	\end{equation}
which shows that the weight operator $W\hk$ can be defined using just the matrices of first $r_k$ singular vectors $\f{U}\hk$ and $\f{V}\hk$ and the $r_k$ first singular values $\sigma_1\hk,\ldots,\sigma_{r_k}\hk$ of $\f{X}\hk$.

For a given rank $r_k$, we recall that the best rank-$r_k$ approximation of a matrix $\f{X}\hk$ as in \cref{eq:Xk:bothsvds} can be written such that
\begin{equation} \label{eq:best:rankr:approximation}
\mathcal{T}_{r_k}(\f{X}\hk) :=  \argmin_{\f{Z}: \rank(\f{Z}) \leq r_k} \|\f{Z} - \f{X}\hk\| = \f{U}\hk \f{\Sigma}\hk \f{V}^{(k)*},
\end{equation}
where $\| \cdot \|$ can be any unitarily invariant norm, due to the Eckardt-Young-Mirsky theorem \cite{Mirsky60}. Let now
\[
\begin{split} 
T_k := T_{\mathcal{T}_{r_k}(\f{X}\hk)}\mathcal{M}_{r_k}: = \{\f{U}\hk \Gamma_1 \f{V}^{(k)*} &+ \f{U}\hk \Gamma_2  \left(\f{I}-\f{V}\hk \f{V}^{(k)*}\right) + (\f{I}-\f{U}^{(k)}\f{U}^{(k)*}) \Gamma_3 \f{V}^{(k)*}: \\ &\Gamma_1 \in \R^{r_k \times r_k}, \Gamma_2 \in \R^{r_k \times d_2}, \Gamma_2\f{V}\hk= 0, \Gamma_3 \in \R^{d_1 \times r_k},  \f{U}^{(k)*}\Gamma_3 = 0  \}
\end{split}
\]
be the tangent space of the manifold of rank-$r_k$ matrices $\mathcal{M}_{r_k}$ of dimension $(d_1 \times d_2)$ at $\mathcal{T}_{r_k}(\f{X}\hk)$ (see \cite{Vandereycken13}), which is a subspace of $\Rdd$ of dimension $r_k(d_1+d_2-r_k)$. Defining 
\begin{equation} \label{def:Sk}
\begin{split}
S_k:= \bigg\{\gamma = (\gamma_1^T,\gamma_2^T,\gamma_3^T)^T \in \R^{r(d_1+d_2+r)}: \Gamma_1=&(\gamma_1)_{\mat} \in \R^{r_k \times r_k}, \Gamma_2=(\gamma_2)_{\mat} \in \R^{r_k \times d_2}, \Gamma_2\f{V}\hk= 0, \\
&\Gamma_3 = (\gamma_3)_{\mat} \in \R^{d_1 \times r_k}, \f{U}^{(k)*}\Gamma_3 = 0\bigg\} \subset \R^{r_k(d_1+d_2+r_k)},
\end{split}
\end{equation}
where $\text{mat}$ is the matricization operator of appropriate dimension that stacks column after column according to the desired dimensions, we can now identify a structure in $W\hk$ that enables us to write it more compactly: Let $P_{T_k}: S_k \to T_k$ be the parametrization operator such that
\[
P_{T_k}(\gamma) := \f{U}\hk \Gamma_1 \f{V}^{(k)*} + \f{U}\hk \Gamma_2  \left(\f{I}-\f{V}\hk \f{V}^{(k)*}\right) + (\f{I}-\f{U}^{(k)}\f{U}^{(k)*}) \Gamma_3 \f{V}^{(k)*}
\]
for $\gamma \in S_k$. Its adjoint operator $P_{T_k}^*: T_k \to S_k$ can be then written such that
\[
P_{T_k}^*(\f{Z}) = \left((\f{Z}_1)_{\vecc}^T,(\f{Z}_2)_{\vecc}^T,(\f{Z}_3)_{\vecc}^T \right)^T
\]
for a matrix $\f{Z}$ such that $T_k \ni \f{Z} =  \f{U}\hk \f{Z}_1 \f{V}^{(k)*} + \f{U}\hk \f{Z}_2  \left(\f{I}-\f{V}\hk \f{V}^{(k)*}\right) + (\f{I}-\f{U}^{(k)}\f{U}^{(k)*}) \f{Z}_3 \f{V}^{(k)*}$, where $\vecc$ is the vectorization operator that concatenates the columns of a matrix, and can be extended to the entire space $\Rdd$ such that
\[
P_{T_k}^*(\f{Z}) = \left((\f{U}_k^* \f{Z} \f{V})_{\vecc}^T,(\f{U}_k^* \f{Z}(\f{I}-\f{V}\hk\f{V}^{(k)*}))_{\vecc}^T,((\f{I}-\f{U}\hk\f{U}^{(k)*})\f{Z} \f{V}\hk)_{\vecc}^T \right)^T
\]
for $\f{Z} \in \Rdd$. With these notations, based on \cref{eq:W:simplification}, we can write $W\hk$ such that
\[
 W\hk = P_{T_k} \left(\f{D}_{S_k}- \epsilon_k^{-2} \f{I}_{S_k}\right)P_{T_k}^* + \epsilon_k^{-2} \f{I},
 \]
 where $\f{I}_{S_k}$ is the identity matrix on $S_k$, $\f{I}$ is the identity operator on $\Rdd$, and $\f{D}_{S_k} \in \R^{r_k(d_1+d_2+r_k) \times r_k(d_1+d_2+r_k)}$ is a diagonal matrix with diagonal entries that are equal to entries of $\f{H}\hk$ from \cref{eq:H:def:small} or to diagonal entries of $\f{D}\hk$ from \cref{eq:H:def:small}, cf. \cref{eq:W:simplification}.
 
 With this preparation, we can now obtain an efficient implementation for the solution of the weighted least squares problem of \cref{eq:MatrixIRLS:Xdef}.
Recall that the set of indices corresponding to provided entries is defined as $\Omega = \{(i_{\ell},j_{\ell})\} \subset [d_1] \times [d_2]$,
and $P_{\Omega}: \Rdd \to \R^m$ is the subsampling operator 
\[
P_{\Omega}(\f{Z}) = \sum_{(i_{\ell},j_{\ell}) \in \Omega} \langle e_{i_{\ell}}, \f{Z} e_{j_{\ell}}\rangle,
\]
where $e_{i_{\ell}}$ and $e_{j_{\ell}}$ are the $i_{\ell}$-th and $j_{\ell}$-th standard basis vectors of $\R^{d_1}$ and $\R^{d_2}$, respectively.

Let $\f{y} \in \R^m$ be the vector of observed entries. As the weight operator $W\hk$ is positive definite, the minimizer 
\[
\f{X}\hkk = \argmin\limits_{P_{\Omega}(\f{X})=\f{y}} \langle \f{X}, W^{(k)}(\f{X}) \rangle
\] 
is unique, and it is well-known \cite{Bjoerck96} that the solution of this linearly constrained weighted least squares problem can be written such that
\begin{equation} \label{eq:weightedleastsquares:formula}
\f{X}\hkk =  (W\hk)^{-1}P_{\Omega}^*\left(P_{\Omega}(W\hk)^{-1} P_{\Omega}^*\right)^{-1}(\f{y}),
\end{equation}
where $(W\hk)^{-1}: \Rdd \to \Rdd$ is the \emph{inverse} of the weight operator $W\hk$: This inverse exists as $W\hk$ is self-adjoint and positive definite, which can be shown by realizing that its $(d_1 d_2 \times d_1 d_2)$-matrix representation has eigenvectors $v_i\hk \otimes u_j\hk$, where $v_i\hk \in \R^{d_1}$ and $u_{j}\hk \in \R^{d_2}$ are columns of $\f{V}\hk$ and $\f{U}\hk$, respectively, and eigenvalues that are just the entries of $\f{H}\hk$.

Furthermore, by replicating the arguments that let to the representation \cref{eq:W:simplification} for the action of $W\hk$ for
\[
(W\hk)^{-1}(\f{Z}) = \f{U}_k \left[\f{H}_k^{-1} \circ (\f{U}_k^{*} \f{Z} \f{V}_k)\right] \f{V}_k^{*},
\]
where $\f{H}_k^{-1} \in \R^{d_1 \times d_2}$ is the matrix whose entries are the \emph{entrywise} inverse of the entries of $\f{H}_k$, it can be seen that we can rewrite $(W\hk)^{-1}$ such that
\begin{equation*}
(W\hk)^{-1} = P_{T_k} \left(\f{D}_{S_k}^{-1}- \epsilon_k^{2} \f{I}_{S_k}\right)P_{T_k}^* + \epsilon_k^{2} \f{I}.
\end{equation*}
Let now $\f{z} \in \R^m$ be the solution of 
\[
\left(P_{\Omega}(W\hk)^{-1} P_{\Omega}^*\right) \f{z} =\f{y}.
\]
Since $P_{\Omega} P_{\Omega}^* = \f{I}_m$, we can use this representation of $(W\hk)^{-1}$ to write
\[
\left(P_{\Omega}(W\hk)^{-1} P_{\Omega}^*\right)^{-1} = \left( P_{\Omega} P_{T_k} \left(\f{D}_{S_k}^{-1}- \epsilon_k^{2} \f{I}_{S_k}\right)P_{T_k}^* P_{\Omega}^* + \epsilon_k^{2} \f{I}_m \right)^{-1}.
\]
Using the \emph{Sherman-Morrison-Woodbury} formula \cite{Woodbury50,Fornasier11} 
\begin{equation*} 
(\f{E} \f{C} \f{F}^* + \f{B})^{-1} = \f{B}^{-1} - \f{B}^{-1} \f{E} ( \f{C}^{-1} + \f{F}^* \f{B}^{-1} \f{E} )^{-1} \f{F}^* \f{B}^{-1}
\end{equation*}
for $\f{B}:= \epsilon_k^{2-p} \f{I}_m$, $\f{C} := \left(\f{D}_{S_k}^{-1}- \epsilon_k^{2} \f{I}_{S_k}\right)$ and $\f{E}= \f{F} := P_{\Omega} P_{T_k}$, we obtain that
\begin{equation} \label{eq:MatrixIRLS:Woodbury:application}
\begin{split}
\f{z} &= \left(P_{\Omega}(W\hk)^{-1} P_{\Omega}^*\right)^{-1}(\f{y}) \\
&= \left[\epsilon_k^{-2} \f{I} - \epsilon_k^{-2}P_{\Omega} P_{T_k}\left(\epsilon_k^{2} \f{C}^{-1} + P_{T_k}^* P_{\Omega}^* P_{\Omega} P_{T_k}\right)^{-1} P_{T_k}^* P_{\Omega}^*\right](\f{y}),
\end{split}
\end{equation}
noting that $\f{C} = \left(\f{D}_{S_k}^{-1}- \epsilon_k^{2} \f{I}_{S_k}\right)$ is invertible as $(\f{H}_{ij}\hk)^{-1} = \sigma_i^{(k)} \sigma_j^{(k)} > \epsilon_k^2$ for all $i,j \in [r_k]$ and as $(\f{D}_{ii}\hk)^{-1} = \sigma_i^{(k)} \epsilon_k > \epsilon_k^2$ for all $i \in [r_k]$, as $r_k = r = \widetilde{r}$ due to the assumption of \Cref{thm:MatrixIRLS:computationalcost:Xkk}.

Inserting this into \cref{eq:weightedleastsquares:formula}, since $P_{T_k}^* P_{T_k} = \f{I}_{S_k}$, and defining the projection operator $\mathcal{P}_{T_k}:= P_{T_k} P_{T_k}^*$ that projects onto $T_k$, we see that
\[
\begin{split}
\f{X}_{T_k}\hkk:= \mathcal{P}_{T_k}(\f{X}\hkk) &:= P_{T_k} P_{T_k}^*(\f{X}\hkk) =  P_{T_k}P_{T_k}^*(W\hk)^{-1}P_{\Omega}^*(\f{z}) =  P_{T_k}\f{D}_{S_k}^{-1} P_{T_k}^*P_{\Omega}^*(\f{z}) \\
&= \epsilon_k^{-2} P_{T_k}\f{D}_{S_k}^{-1} \left[P_{T_k}^*P_{\Omega}^* - P_{T_k}^*P_{\Omega}^* P_{\Omega} P_{T_k}\left(\epsilon_k^{2} \f{C}^{-1} + P_{T_k}^* P_{\Omega}^* P_{\Omega} P_{T_k}\right)^{-1} P_{T_k}^* P_{\Omega}^*\right](\f{y}) \\
&= \epsilon_k^{-2} P_{T_k}\f{D}_{S_k}^{-1} \left[\epsilon_k^{2} \f{C}^{-1}\left(\epsilon_k^{2} \f{C}^{-1} + P_{T_k}^* P_{\Omega}^* P_{\Omega} P_{T_k}\right)^{-1} P_{T_k}^* P_{\Omega}^*\right](\f{y}) \\
&= P_{T_k}\f{D}_{S_k}^{-1} \left(\f{D}_{S_k}^{-1}- \epsilon_k^{2} \f{I}_{S_k}\right)^{-1} \gamma_k,
\end{split}
\]
where $\gamma_k \in S_k$ is the solution of the linear system of size $(r_k (d_1 + d_2 + r_k) \times r_k (d_1 + d_2 +r_k))$
\begin{equation} \label{eq:gamma:system}
\left(\epsilon_k^{2} \left(\f{D}_{S_k}^{-1}- \epsilon_k^{2} \f{I}_{S_k}\right)^{-1} + P_{T_k}^* P_{\Omega}^* P_{\Omega}  P_{T_k}\right) \f{\gamma}_{k} =  P_{T_k}^* P_{\Omega}^* (\f{y}).
\end{equation}
For the part $\f{X}_{T_k^{\perp}}\hkk:= \left(\f{I}- \mathcal{P}_{T_k} \right)\f{X}\hkk$ of $\f{X}\hkk$ in the orthogonal complement $T_k^\perp$ of $T_k$, we observe that
\[
\begin{split}
\f{X}_{T_k^{\perp}}\hkk &= \left(\f{I}- P_{T_k} P_{T_k}^* \right)(W\hk)^{-1}P_{\Omega}^*(\f{z})= \epsilon_k^2 \left(\f{I}- P_{T_k} P_{T_k}^* \right)  P_{\Omega}^*(\f{z})  \\
&= \left(\f{I}- P_{T_k} P_{T_k}^* \right)P_{\Omega}^*  \left[\f{I} - P_{\Omega} P_{T_k}\left(\epsilon_k^{2} \f{C}^{-1} + P_{T_k}^* P_{\Omega}^* P_{\Omega} P_{T_k}\right)^{-1} P_{T_k}^* P_{\Omega}^*\right](\f{y}) = \\
&= \left(\f{I}- P_{T_k} P_{T_k}^* \right)P_{\Omega}^*  \left(\f{y} - P_{\Omega} P_{T_k}(\f{\gamma}_{k}) \right) = P_{\Omega}^* (\f{r}_{k+1}) - P_{T_k} P_{T_k}^*P_{\Omega}^*(\f{r}_{k+1})
\end{split}
\]
with the definition of $\f{r}_{k+1} := \f{y} - P_{\Omega} P_{T_k}(\f{\gamma}_{k}) \in \R^m$.

Combining these two parts of $\f{X}\hkk$, we see that the solution $\f{X}\hkk$ of the linearly constrained weighted least squares problem can be written as
\[
\begin{split}
\f{X}\hkk &= P_{\Omega}^* (\f{r}_{k+1}) + P_{T_k} \left(\left(\f{D}_{S_k}^{-1}- \epsilon_k^{2} \f{I}_{S_k}\right)^{-1} \gamma_k -P_{T_k}^*P_{\Omega}^*(\f{r}_{k+1}) \right) =: P_{\Omega}^* (\f{r}_{k+1}) + P_{T_k} (\widetilde{\gamma}_k) \\
&= P_{\Omega}^* (\f{r}_{k+1}) + \f{U}\hk \left(\widetilde{\Gamma}_1 \f{V}^{(k)*} + \widetilde{\Gamma}_2\right) + \widetilde{\Gamma}_3 \f{V}^{(k)*} \\
&= P_{\Omega}^* (\f{r}_{k+1}) + \f{U}\hk \f{M}_{1}^{(k+1)*} + \f{M}_2^{(k+1)} \f{V}^{(k)*},
\end{split}
\]
which provides a representation of $\f{X}\hkk$ as a sum of the $m$-sparse matrix $P_{\Omega}^* (\f{r}_{k+1})$ and the two matrices $\f{U}\hk \f{M}_{1}^{(k+1)*}$ and $\f{M}_2^{(k+1)} \f{V}^{(k)*}$  of rank $r=r_k$, as stated in \Cref{thm:MatrixIRLS:computationalcost:Xkk}, if we define $\widetilde{\gamma}_k \in S_k$ as in the first line, defining $\widetilde{\Gamma}_1 \in \R^{r_k \times r_k}$, $\widetilde{\Gamma}_2 \in \R^{r_k \times d_2}$ and $\widetilde{\Gamma}_3 \in \R^{d_1 \times r_k}$ as the representation matrices of $\widetilde{\gamma}_k$ in \cref{def:Sk}, and 
\[
\f{M}_{1}^{(k+1)}:= \f{V}\hk \widetilde{\Gamma}_1^* + \widetilde{\Gamma}_2^*, \quad \f{M}_{2}^{(k+1)}:= \widetilde{\Gamma}_3.
\]
Evidently, this representation of $\f{X}\hkk$ has $m + r_k d_2 +d_1 r_k = O(m+ rD)$ parameters.
 
As a summary, we implement the solution of the weighted least squares problem \cref{eq:weightedleastsquares:formula} as follows:
\begin{enumerate}
	\item Calculate $P_{T_k}^*P_{\Omega}^*(\f{y}) \in S_k$.
	\item Solve \cref{eq:gamma:system} for $\gamma_k \in S_k$ by the \emph{conjugate gradient} method \cite{HestenesStiefel52,Meurant}.
	\item Calculate residual $\f{r}_k:= \f{y}-P_{\Omega} P_{T_k} (\f{\gamma}_{k}) \in \R^m$.
	\item Calculate $\widetilde{\f{\gamma}_k} = \left(\f{D}_{S_k}^{-1}- \epsilon_k^{2} \f{I}_{S_k}\right)^{-1} \gamma_k - P_{T_k}^* P_{\Omega}^* (\f{r}_k) \in S_k$ and set $\f{M}_{1}^{(k+1)}:= \f{V}\hk \widetilde{\Gamma}_1^* + \widetilde{\Gamma}_2^*$ and $\f{M}_{2}^{(k+1)}:= \widetilde{\Gamma}_3$.
\end{enumerate}
The main computational cost lies in the application of the operators $(P_{\Omega} P_{T_k})^* = P_{T_k}^* P_{\Omega}^*:\R^m \to S_k$, $P_{\Omega} P_{T_k}: S_k \to \R^m$ and $\left(\f{D}_{S_k}^{-1}- \epsilon_k^{2} \f{I}_{S_k}\right)^{-1} : S_k \to S_k$. The application of the first has a time complexity of $2 r m + r d+ r^2 (d_1 + d_2) = O(m r + r^2 D)$, and $P_{\Omega} P_{T_k}(\gamma)$ can be calculated such that for each $\ell \in [m]$,
\[
\left(P_{\Omega} P_{T_k}(\gamma)\right)_{\ell} = \left(P_{\Omega}\left( \f{U}\hk (\Gamma_1 \f{V}^*+\Gamma_2) + \Gamma_3 \f{V}^{(k)*}\right)\right)_{\ell} = \sum_{k=1}^r (\f{U}\hk)_{i_{\ell},k} \left((\Gamma_1 \f{V}^{(k)*})_{k,j_{\ell}}+(\Gamma_2)_{k,j_{\ell}}\right) + \sum_{k=1}^r (\Gamma_3)_{i_{\ell},k} (\f{V}^{(k)*})_{k,j_{\ell}},
\]
resulting in a number of $O(m r + r^2 D)$ flops. Furthermore, $\left(\f{D}_{S_k}^{-1}- \epsilon_k^{2} \f{I}_{S_k}\right)^{-1}$ is diagonal, requiring $O(r D)$ flops if applied to an element of $S_k$.

In particular, an iterative solver of \cref{eq:gamma:system} such as the conjugate gradient (CG) method, which is applicable since its system matrix positive definite, applies these three operators at each inner iteration. It is known that in general, the CG method terminates with the exact solution $\gamma_k$ after at most $N_{\text{CG\_inner}} = \dim(S_k)= O(r D)$ iterations. 

However, the CG method returns very high precision approximate solutions after a constant number of iterations $N_{\text{CG\_inner}}$ if the system matrix 
\[
\frac{\epsilon_k^{2} \f{I}}{\f{D}_{S_k}^{-1}- \epsilon_k^{2} \f{I}_{S_k}} + P_{T_k}^* P_{\Omega}^* P_{\Omega}  P_{T_k} 
\]
is well-conditioned. And indeed, it can be verified that the spectrum of $\frac{\epsilon_k^{2} \f{I}}{\f{D}_{S_k}^{-1}- \epsilon_k^{2} \f{I}_{S_k}}$ converges to $0$ as $\epsilon_k \to 0$, which is the case as $\f{X}\hk$ converges to a rank-$r$ matrix by virtue of the choice of the smoothing \cref{eq:MatrixIRLS:epsdef}. We note that it can be shown that the spectrum of $P_{T_k}^* P_{\Omega}^* P_{\Omega}  P_{T_k}$ concentrates around $\frac{m}{d_1 d_2}$ if $\mathcal{T}_{r_k}(\f{X}\hk)$ is $\mu_0$-incoherent and if $m = \Omega( \mu_0 r (d_1 + d_2) \log(d_1 + d_2))$ under a typical random measurement model on the indices $\Omega$ (see, e.g., Theorem 6 of \cite{recht}), implying that a constant condition number of the system matrix of \cref{eq:gamma:system} can be shown in this case. 

Overall, by reformulating the weighted least squares problem such that the linear system to solve is \cref{eq:gamma:system}, we avoid a bad conditioning of the system matrix to be solved, which is a common problem of IRLS methods \cite{Fornasier16}, as a naive implementation would suffer from a \emph{blow up} of the weights (in our case, the weight operator $W\hk$) as $\epsilon_k \to 0$.

\section{Remarks to MatrixIRLS as a saddle-escaping smoothing Newton method} \label{sec:remarks:smoothing:newton}
We briefly elaborate on the interpretation of \texttt{MatrixIRLS} as a saddle-escaping smoothing method.

If $\epsilon_k > 0$ and if $F_{\epsilon_k}:\Rdd \to R$ is the $\epsilon_k$-smoothed log-det objective of \cref{eq:smoothing:Fpeps}, it can be shown that $F_{\epsilon_k}$ is continuously differentiable with $\epsilon_k^{-2}$-Lipschitz gradient 
\begin{equation*} 
 \nabla F_{\epsilon_k}(\f{X}) = \f{U} \dg \bigg(\frac{\sigma_i(\f{X})}{\max(\sigma_i(\f{X}),\epsilon_k)^{2}}\bigg)_{i=1}^d \f{V}^*
\end{equation*}
for any matrix $\f{X}$ with singular value decomposition $\f{X} = \f{U} \dg\big(\sigma(\f{X})\big) \f{V}^* = \f{U} \dg\big(\sigma \big) \f{V}^*$. This can be shown by using results from \cite{Lewis05_Nonsm1,AnderssonCarlssonPerfekt16}. Additionally, it is holds that $\nabla F_{\epsilon_k}$ is differentiable at $\f{X}$ if and only if the derivative $f_{\epsilon_k}':\R \to \R$ of $f_{\epsilon_k}$ from \cref{eq:smoothing:Fpeps} exists at all $\sigma=\sigma_i(\f{X})$, $i \in [d]$, which is the case if $\f{X} \in \mathcal{D}_{\epsilon_k}:=\big\{ \f{X}: \sigma_i (\f{X}) \neq \epsilon_k \text{ for all } i \in [d] \big\}$. The latter statement follows from the calculus of \emph{non-Hermitian L\"owner functions} \cite{Yang09,Ding18}, also called \emph{generalized matrix functions} \cite{Noferini17}, as $\f{X} \mapsto \nabla F_{\epsilon_k}(\f{X})$ is such a function.

Let now $\f{X}\hk \in \mathcal{D}_{\epsilon_k}:=\big\{ \f{X}: \sigma_i (\f{X}) \neq \epsilon_k \text{ for all } i \in [d] \big\}$ with singular value decomposition as in \cref{eq:Xk:bothsvds}, and $r_k := |\{i \in [d]: \sigma_i(\f{X}\hk) > \epsilon_k\}| = |\{i \in [d]: \sigma_i\hk > \epsilon_k\}|$. In this case, it can be calculated that the Hessian $\nabla^2 F_{\epsilon_k}(\f{X}\hk)$ at $\f{X}\hk$, which is a function that maps $\Rdd$ to $\Rdd$ matrices, satisfies in the case of $d_1=d_2$
\begin{equation} \label{eq:smoothedranksurrogate:Hessianformula}
\nabla^2 F_{\epsilon_k}(\f{X}\hk)(\f{Z}) =  \f{U}_k \begin{bmatrix}\f{M}^{\text{S}} \circ S(\f{U}_k^* \f{Z} \f{V}_{k})+ 	
	\f{M}^{\text{T}} \circ T(\f{U}_k^* \f{Z} \f{V}_k) \end{bmatrix}
		 \f{V}_k^*, 
\end{equation}
for any $\f{Z} \in \Rdd$, where $S:\R^{d \times d} \to \R^{d \times d}$ is the \emph{symmetrization operator} that maps any $\f{X} \in \R^{d \times d}$ such that
\begin{equation} \label{eq:symmetrization}
	S(\f{X})=\frac{1}{2}(\f{X} + \f{X}^*),
\end{equation}
and  $T:\R^{d \times d} \to \R^{d \times d}$ is the \emph{antisymmetrization operator} such that
\begin{equation} \label{eq:antisymmetrization}
 T(\f{X})= \frac{1}{2}(\f{X} - \f{X}^*)
\end{equation}
for any $\f{X} \in \R^{d \times d}$, and  $\f{M}^{\text{S}}, \f{M}^{\text{T}}  \in \Rdd$ fulfill 
\[
\f{M}^{\text{S}}=
\left[
\begin{array}{c|c}
-\f{H}\hk & \f{M}_{1,2}^{-} \\
\hline
\f{M}_{2,1}^{-} & \epsilon_k^{-2} \mathbf{1}
\end{array}
\right] \quad \quad 
\f{M}^{\text{T}}=
\left[
\begin{array}{c|c}
-\f{H}\hk & \f{M}_{1,2}^{+} \\
\hline
\f{M}_{2,1}^{+} & \epsilon_k^{-2} \mathbf{1}
\end{array}
\right]
\]
with $\f{H}\hk \in \R^{r_k \times r_k}$ as in \cref{eq:H:def:small} and the $(d_1-r_k) \times (d_2- r_k)$-matrix of ones $\mathbf{1}$. Furthermore, the matrices $\f{M}_{1,2}^{-},\f{M}_{1,2}^{+} \in (d_1- r_k) \times r_k$ are such that
\[
\left(\f{M}_{1,2}^{\pm}\right)_{ij} = \frac{(\sigma_i^{(k)})^{-1} \pm \sigma_{j+r_k}^{(k)}\epsilon_k^{-2}}{\sigma_i^{(k)} \pm \sigma_{j+r_k}^{(k)}} 
\] 
for $i \in [r_k]$, $j \in [d_2-r_k]$ and
\[
\left(\f{M}_{2,1}^{\pm}\right)_{ij} = \frac{(\sigma_j^{(k)})^{-1} \pm \sigma_{i+r_k}^{(k)}\epsilon_k^{-2}}{\sigma_j^{(k)} \pm \sigma_{i+r_k}^{(k)}} 
\] 
for$j \in [r_k]$, $i \in [d_1-r_k]$. The formula \cref{eq:smoothedranksurrogate:Hessianformula} for $\nabla^2 F_{\epsilon_k}(\f{X}\hk)$ follows by inserting into Theorem 2.2.6 of \cite{Yang09}, Corollary 3.10 \cite{Noferini17} or Theorem 4 of \cite{Ding18}.

By realizing that $0 \leq \sigma_{\ell}\hk \leq \epsilon_k$ for all $\ell > r_k$, we see that
\[
\frac{1}{(\sigma_i\hk)^{2}} \leq \left(\f{M}_{1,2}^{+}\right)_{ij} = \left(\f{M}_{2,1}^{+}\right)_{ji} \leq  \frac{1}{\sigma_i\hk \epsilon_k}
\]
and 
\[
-\frac{1}{\sigma_i\hk \epsilon_k} \leq \left(\f{M}_{1,2}^{-}\right)_{ij} = \left(\f{M}_{2,1}^{-}\right)_{ji} \leq  \frac{1}{(\sigma_i\hk)^{2}}
\]
for all $i$ and $j$.

Now, comparing $M^{\text{S}}$ and $M^{\text{T}}$ with $\f{H}_k$, see \cref{eq:Hk:blockstructure}, of the weight operator $W\hk$, we see that the upper left blocks of $M^{\text{S}}$ and $M^{\text{T}}$ are just the \emph{negative} of the upper left block $\f{H}\hk$ of $\f{H}_k$, while the lower right blocks coincide. Furthermore, the lower left and the upper right blocks are related such that
\[
\left|\left(\f{M}_{1,2}^{\pm}\right)_{ij}\right| \leq \frac{1}{\sigma_i\hk \epsilon_k} = (\f{H}_{1,2}\hk)_{ij}
\]
for all $i \in [r_k]$, $j \in [d_2-r_k]$, and 
\[
\left|\left(\f{M}_{2,1}^{\pm}\right)_{ij}\right| \leq  \frac{1}{\sigma_j\hk \epsilon_k} = (\f{H}_{2,1}\hk)_{ij}
\]
for all $i \in [d_2-r_k]$, $j \in [r_k]$.

We now point out the relationship of these considerations to an analysis that was performed in \cite{PaternainMokhtariRibeiro19} for the case of an \emph{unconstrained minimization} of $F_{\epsilon_k}$, assuming furthermore that $F_{\epsilon_k}$ was smooth:

In this case, \cite{PaternainMokhtariRibeiro19} considers using \emph{modified Newton steps} 
\[
\f{X}\hkk := \f{X}\hk - \eta_k \left|\nabla^2 F_{\epsilon_k}(\f{X}\hk)\right|_{c}^{-1} \nabla F_{\epsilon_k}(\f{X}\hk)
\]
where the Hessian $\nabla^2 F_{\epsilon_k}(\f{X}\hk)$ is replaced by a positive definite truncated eigenvalue matrix $\left|\nabla^2 F_{\epsilon_k}(\f{X}\hk)\right|_{c}$, which replaces the large negative eigenvalues of $\nabla^2 F_{\epsilon_k}(\f{X}\hk)$ by their modulus for eigenvalues that have large modulus and eigenvalues of small modulus by an appropriate constant $c$. \cite{PaternainMokhtariRibeiro19} shows that such steps are, unlike conventional Newton steps, which often are \emph{attracted by saddle points}, able to \emph{escape} saddle points with an exponential rate that does \emph{not} depend on the conditioning of the problem. Experimental observations of such behavior has been reported also in other works \cite{Murray2010,Dauphin14}.

In view of this, we observe that the weight operator $W\hk$ is nothing but a refined variant of $\left|\nabla^2 F_{\epsilon_k}(\f{X}\hk)\right|_{c}$, as the eigenvalues of $\nabla^2 F_{\epsilon_k}(\f{X}\hk)$ from \cref{eq:smoothedranksurrogate:Hessianformula} are simply $\{ (\f{M}^{\text{S}}_{ij}, i \leq j\} \cup \{ (\f{M}^{\text{T}}_{ij}, i < j\}$, c.f., e.g., Theorem 4.5 of \cite{Noferini17}. In particular, the refinement is such that the small eigenvalues of $\nabla^2 F_{\epsilon_k}(\f{X}\hk)$, which can be found in the entries of $\f{M}_{1,2}^{\pm}$ and $\f{M}_{2,1}^{\pm}$, are replaced not by a uniform constant, but by \emph{different} upper bounds $(\sigma_i\hk \epsilon_k)^{-1}$ and $(\sigma_j\hk \epsilon_k)^{-1}$ that depend either on the row index $i$ or the column index $j$.

Besides this connection, there are important differences of our algorithm to the algorithm analyzed in \cite{PaternainMokhtariRibeiro19}. While that paper considers the minimization of a fixed smooth function, we update the smoothing parameter $\epsilon_k$ and thus the function $F_{\epsilon_k}$ at each iteration. Furthermore, Algorithm 1 of \cite{PaternainMokhtariRibeiro19} uses backtracking for each modified Newton step, which would be prohibitive to perform as evaluations of $F_{\epsilon_k}$ are very expensive for our smoothed log-det objectives, as they would require the calculation of all singular values. On the other hand, \texttt{MatrixIRLS} uses full modified Newton steps, and we can assure that these are always a decent direction in our case, as we explain in an upcoming paper. Lastly, we do not add noise to the iterates.

As mentioned in \Cref{surrogate}, \texttt{MatrixIRLS} is by no means the \emph{first} algorithm for low-rank matrix recovery that can be considered as an iteratively reweighted least squares algorithm. However, the IRLS algorithms \cite{Fornasier11,Mohan10,Lai13,KS18} are different from \texttt{MatrixIRLS} not only in their computational aspects, but also since they do \emph{not} allow for a close relationship between their weight operator $W\hk$ and the Hessian $\nabla^2 F_{\epsilon_k}(\f{X}\hk)$ at $\f{X}\hk$ as described above.

\section{Experimental Details} \label{sec:experimental:details}
In this section, we specify some details of the setup and the algorithmic parameters of the experiments presented in \Cref{sec:numerics}.

First, we note that if we are interested in recovering a rank-$r$ matrix $\f{X}^0 \in \Rdd$ from its entries indexed by a set $\Omega$ for small oversampling factors $\rho$ and corresponding sample sizes $|\Omega|=: m = \lfloor \rho r (d_1 +d_2 - r) \rfloor$, the solution of even the intractable rank minimization formulation
\begin{equation} \label{eq:rank:min:POmegaX0}
\min_{\f{X} \in \Rdd} \rank(\f{X})  \quad \mbox{ subject to } P_{\Omega}(\f{X}) = P_{\Omega}(\f{X}^0)
\end{equation}
might \emph{not} coincide with $\f{X}^0$, or the solution might not be unique, even if the sample set $\Omega$ is chosen uniformly at random.

In particular, this will be the case if $\Omega$ is such that there is a row or a column with \emph{fewer than} $r$ revealed entries, as a necessary condition for the unique reconstruction by \cref{eq:rank:min:POmegaX0} would be violated in this case \cite{Pimentel15}.

To mitigate this problem that is rather related to the structure of the sampling set than to the performance of a certain algorithm, we, in fact, adapt the sampling model of uniform sampling without replacement by checking whether the condition such that each row and each column in $\Omega$ has at least $r$ observed entries, and resampling $\Omega$ if this condition is not fulfilled. This procedure is repeated up to a maximum of $1000$ resamplings.

The sample complexity experiments were conducted on a Linux node with Intel Xeon E5-2690 v3 CPU with 28 cores and 64 GB RAM, using MATLAB R2019a. All the other experiments were conducted on a Windows 10 laptop with Intel i7 7660U with 2 cores and 8 GB RAM, also using MATLAB R2019a. For the purpose of our experiment, we categorize the algorithms into algorithms of \emph{first-order} type and of \emph{second-order} type based on whether an algorithm exhibits empiricially observed locally superlinear convergence rates or not.

 All the methods are provided with the true rank $r$ of $\f{X}^0$ as an input parameter. If possible, we use the MATLAB implementation provided by the authors of the respective papers. We do not make use of explicit parallelization for any of the methods, but most methods use complied C subroutines to efficiently implement sparse evaluations of matrix factorizations.

We set a maximal number of outer iterations for the second-order methods as $N_0 = 400$. The second-order type algorithms we consider are:
\begin{itemize}
	\item \texttt{MatrixIRLS}, as described in \Cref{algo:MatrixIRLS}. As a stopping criterion, we choose a threshold of $10^{-9}$ for the relative change of the Frobenius norm $\frac{\|\f{X}\hkk - \f{X}\hk\|_F}{\|\f{X}\hk\|_F}$. We use the CG method for solving the linear system \cref{eq:gamma:system} without any preconditioning, but using 
	\[
	\gamma_{k,0}:= \left(\f{D}_{S_k}^{-1}- \epsilon_k^{2} \f{I}_{S_k}\right) P_{T_k}^*(\f{X}\hk)
	\]
	as an initial guess for the solution of \cref{eq:gamma:system}. We terminate the CG method if a maximum number of $N_{\text{CG\_{inner}}} = 500$ inner iterations is reached or if a relative residual of $tol_{\text{inner}} = 10^{-5} \epsilon_{k}\kappa^{-1}$ is reached, whichever happens first.\footnote{While this stopping condition uses the condition number $\kappa$, which will probably be unknown in practice, it can be generally chosen independently of $\kappa$ without any problems of convergence.} For the weight operator update step, we use a variant of the randomized Block Krylov method \cite{MuscoMusco15} based on the implementation provided by the authors\footnote{\url{https://github.com/cpmusco/bksvd}}, setting the parameter for the maximal number of iterations to $20$.
	\item \texttt{R2RILS} \cite{BauchNadler20} or \emph{rank $2r$ iterative least squares}, a method that optimizes a least squares data fit objective $\|P_{\Omega}(\f{X}_{0}) - P_{\Omega}(\f{X})\|_F$ over $\f{X} \in T_{\f{Z}\hk}\mathcal{M}_{r}$, where $T_{\f{Z}\hk}\mathcal{M}_{r}$ is a tangent space onto the manifold of rank-$r$ matrix manifold, while iteratively updating this tangent space. As above, we stop the outer iterations a threshold of $10^{-9}$ is reached for the relative change of the Frobenius norm $\frac{\|\f{X}\hkk - \f{X}\hk\|_F}{\|\f{X}\hk\|_F}$. At each outer iteration, \texttt{R2RILS} solves an overdetermined least squares problem of size $(m \times r(d_1 +d_2))$ via the iterative solver LSQR, for which we choose the maximal number of inner iterations as $N_{\text{LSQR\_{inner}}} = 500$ and a termination criterion based on a relative residual of $10^{-6} \kappa^{-1}$. We use the implementation based on the code provided by the authors, but adapted for these stopping criteria.\footnote{\url{https://github.com/Jonathan-WIS/R2RILS}}
	\item  \texttt{RTRMC}, the preconditioned Riemannian trust-region method called RTRMC 2p of \cite{boumal_absil_15}, which was reported to achieve the best performance among a variety of matrix completion algorithms for the task of completing matrices of a condition number of up to $\kappa = 150$. We use the implementation provided by the authors\footnote{RTRMC v3.2 from \url{http://web.math.princeton.edu/~nboumal/RTRMC/index.html}, together with the toolbox Manopt 6.0 (\url{https://www.manopt.org/}) \cite{manopt}.} with default options except from setting the maximal number of inner iterations to $N_{\text{inner}}=500$ and setting the parameter for the tolerance on the gradient norm to $10^{-15}$. Furthermore, as the algorithm otherwise would often run into certain submatrices that are not positive definite for $\rho$ between $1$ and $1.5$, we set the regularization parameter $\lambda = 10^{-8}$, which is small enough not to deter high precision approximations of $\f{X}^0$ if enough samples are provided.
\end{itemize}

Furthermore, we consider the following first-order algorithms, setting the maximal number of outer iterations to $N_0=4000$:

\begin{itemize}
	\item \texttt{LRGeomCG} \cite{Vandereycken13}, a local optimization method for a quadratic data fit term based on gradients with respect to the Riemannian manifold of fixed rank matrices. We use the author's implementation\footnote{\url{http://www.unige.ch/math/vandereycken/matrix_completion.html}} while setting the parameters related to the stopping conditions \texttt{abs\_grad\_tol}, \texttt{rel\_grad\_tol}, \texttt{abs\_f\_tol}, \texttt{rel\_f\_tol}, \texttt{rel\_tol\_change\_x} and \texttt{rel\_tol\_change\_res} each to $10^{-9}$.
	\item  \texttt{LMaFit} or low-rank matrix fitting \cite{Wen12}, a nonlinear successive over-relaxation algorithm based on matrix factorization. We use the implementation provided by the authors\footnote{\url{http://lmafit.blogs.rice.edu}}, setting the tolerance threshold for the stopping condition (which is based on a relative data fit error $\|P_{\Omega}(\f{X}^{(k)})-\f{y}\|_2/\|\f{y}\|_2$) to $5 \cdot 10^{-10}$.
	\item \texttt{ScaledASD} or scaled alternating steepest descent \cite{TannerWei16}, a gradient descent method based on matrix factorization that which scales the gradients in a quasi-Newton fashion. We use the implementation provided by the authors\footnote{\label{algWei}\url{http://www.sdspeople.fudan.edu.cn/weike/code/mc20140528.tar}} with the stopping condition of $\|P_{\Omega}(\f{X}^{(k)})-\f{y}\|_2/\|\f{y}\|_2 \leq 10^{-9}$.
	\item  \texttt{ScaledGD} or scaled gradient descent \cite{tong_ma_chi}, a method that is very similar to \texttt{ScaledASD}, but for which a non-asymptotic local convergence analysis has been achieved for the case of a matrix recovery problem related to matrix completion, and which has been investigated experimentally in \cite{tong_ma_chi} in the light of the completion of ill-conditioned low-rank matrices. We use an adapted version of the author's implementation\footnote{\url{https://github.com/Titan-Tong/ScaledGD}}: We choose a step size of $\eta = 0.5$, but increase the normalization parameter $p$ by a factor of $1.5$ in case the unmodified algorithm \texttt{ScaledGD} leads to divergent algorithmic iterates, using the same stopping condition as for \texttt{ScaledASD}.
	\item \texttt{NIHT} or normalized iterative hard thresholding \cite{TannerWei13}, which performs iterative hard thresholding steps with adaptive step sizes. We use the implementation provided by the authors \footref{algWei} with a stopping threshold of $10^{-9}$ for the relative data fit error $\|P_{\Omega}(\f{X}^{(k)})-\f{y}\|_2/\|\f{y}\|_2$ and the convergence rate threshold parameter $1- 10^{-9}$. 
\end{itemize}

We base our choice of algorithms on the desire to obtain a representative picture of state-of-the-art algorithms for matrix completion, including in particular those that are scalable to problems with dimensionality in the thousands or more, those that come with the best theoretical guarantees, and those that claim to perform particularly well to complete \emph{ill-conditioned} matrices.

\begin{centering}
\begin{figure}[ht]
 \includegraphics[width=0.8\textwidth,height=0.4\textheight,clip]{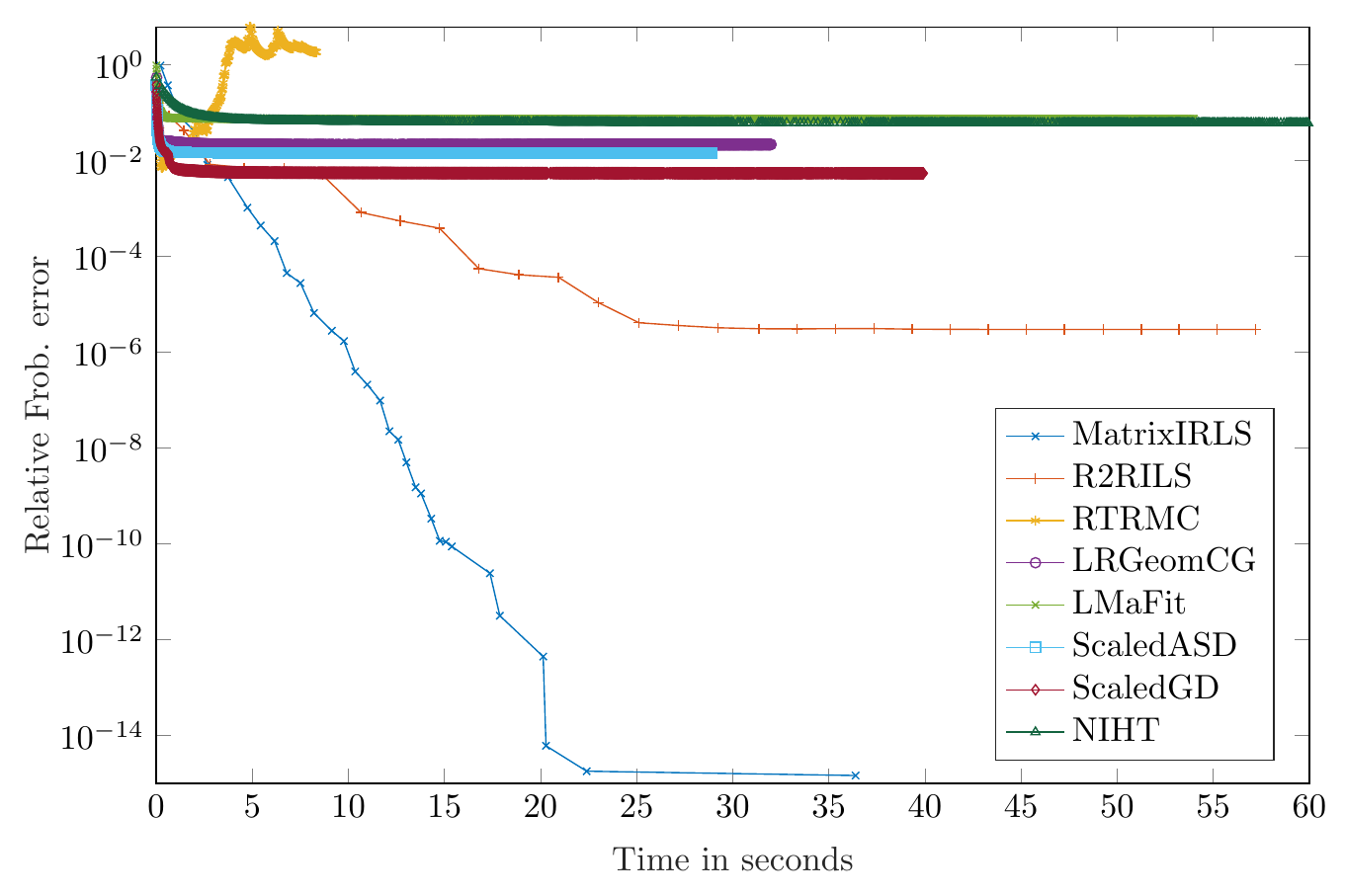}
 \centering
\vspace*{-5mm}
\caption{Completion of a highly ill-conditioned $1000 \times 1000$ matrix of rank $r=10$ with $\kappa=10^{10}$ with an oversampling rate of $\rho= 4$.}
\label{running_time_extremeill_cond}
\end{figure}
\end{centering}

As an additional experiment, in Figure \ref{running_time_extremeill_cond} we performed simulations for extremely ill-conditioned $1000 \times 1000$ matrices with $\rank=10 $ and $\kappa=10^{10}$ and exponentially interpolated singular values as described above. For this case, \texttt{MatrixIRLS} attains a relative Frobenius error of the order of the machine precision and, remarkably, exactly recover all the singular values up to 15 digits. This also shows that the conjugated gradient and the randomized block Krylov method used at the inner core of our implementation can be extremely precise when properly adjusted. \texttt{R2RILS} is also able to obtain relatively low Frobenius error but unlike our method, it is not able to retrieve all the singular values with high accuracy. All the other methods we ran were not able to complete such extremely ill-conditioned matrices.

\vspace{-.5cm}

\begin{table}[h]
\begin{tabular}{|l|l|l|l|l|}
\hline
$\sigma_1=1000000$     & $\sigma_2=77426.368268113$ & $\sigma_3=5994.842503189$ & $\sigma_4=464.158883361$ & $\sigma_5=35.938136638$   \\ \hline
$\sigma_6=2.782559402$ & $\sigma_7=0.215443469$     & $\sigma_8=0.016681005$    & $\sigma_9=0.001291550$   & $\sigma_{10}=0.000100000$ \\ \hline
\end{tabular}
\caption{Singular values of the $1000 \times 1000$ random matrix used in the experiment of Figure \ref{running_time_extremeill_cond}. All the singular values are multiplied by $10^{-4}$ in order to fit the table inside the page. As one can observe, $\sigma_1=10^{10}$ and $\sigma_{10}=1$, while all the other values are exponentially interpolated in between. \texttt{MatrixIRLS} is able to recover all of them with very high accuracy.}
\label{tab:table_extreme_ill}
\end{table}

\end{document}